\nonstopmode
\documentclass[a4paper, 10pt]{amsart}
\usepackage{latexsym}
\usepackage{fancyhdr}
\usepackage{amsmath, amssymb}
\usepackage[ansinew]{inputenc}
\usepackage[all]{xy}
\usepackage{verbatim}
\usepackage{psfrag}
\usepackage{epsfig}
\swapnumbers
\theoremstyle{plain}
\newtheorem*{lemma*}{Lemma}
\newtheorem{lemma}[subsection]{Lemma}
\newtheorem*{theorem*}{Theorem}
\newtheorem{theorem}[subsection]{Theorem}
\newtheorem*{proposition*}{Proposition}
\newtheorem{proposition}[subsection]{Proposition}
\newtheorem*{corollary*}{Corollary}
\newtheorem{corollary}[subsection]{Corollary}
\theoremstyle{definition}
\newtheorem*{definition*}{Definition}

\newtheorem*{example*}{Example}

\newtheorem*{remark*}{Remark}
\newtheorem*{remarks*}{Remarks}

\numberwithin{equation}{subsection}
\newenvironment{demo}[1]{\par\smallskip\noindent{\bf #1.}}{\par\smallskip}

\allowdisplaybreaks

\let\on=\operatorname
\newcommand{\sr}[1]%
{\ifmmode{}^\dagger\else${}^\dagger$\fi\ifvmode
\vbox to 0pt{\vss
\hbox to 0pt{\hskip\hsize\hskip1em
\vbox{\hsize3cm\raggedright\pretolerance10000
\noindent #1\hfill}\hss}\vss}\else
\vadjust{\vbox to0pt{\vss%
 \hbox to 0pt{\hskip\hsize\hskip1em%
 \vbox{\hsize3cm\raggedright\pretolerance10000%
 \noindent #1\hfill}\hss}\vss}}\fi%
}

\ifcase\pdfoutput

\else
\usepackage[pdftex,raiselinks=true,backref=page
]{hyperref}%
\hypersetup{
bookmarksnumbered=true,%
bookmarksopen=true,%
bookmarksopenlevel=1
}
\fi

\providecommand{\mapsfrom}{\kern.2em%
\setbox0=\hbox{$\leftarrow$\kern-.10em\rule[0.26mm]{0.1mm}{1.3mm}}\box0%
\kern.3em}

\allowdisplaybreaks
\date{\today}

\title[Quasianalytic mappings]
{The Convenient Setting for Quasianalytic Denjoy--Carleman
Differentiable Mappings}

\author
[A.~Kriegl, P.W.~Michor, A.~Rainer]
{Andreas Kriegl, Peter W. Michor, and Armin Rainer}

\address{Andreas Kriegl: Fakult\"at f\"ur Mathematik, Universit\"at Wien, Nordbergstrasse~15, A-1090 Wien, Austria}
\email{andreas.kriegl@univie.ac.at}

\address{Peter W. Michor: Fakult\"at f\"ur Mathematik, Universit\"at Wien, Nordbergstrasse~15, A-1090 Wien, Austria}
\email{peter.michor@univie.ac.at}

\address{Armin Rainer: Department of Mathematics, University of Toronto, 40 St.\ George Street, Toronto, Ontario, Canada M5S 2E4}

\email{armin.rainer@univie.ac.at}

\thanks{PM was supported by FWF-Project P~21030-N13.
AR was supported by FWF-Project J2771}

\begin{document}

\begin{abstract}
For   quasianalytic  Denjoy--Carleman  differentiable  function
classes   $C^Q$   where   the   weight  sequence  $Q=(Q_k)$  is
log-convex,  stable  under  derivations, of moderate growth and
also an $\mathcal L$-intersection (see \thetag{\ref{nmb:1.6}}),
we  prove  the  following:  The  category  of $C^Q$-mappings is
cartesian   closed  in  the  sense  that  $C^Q(E,C^Q(F,G))\cong
C^Q(E\times  F,  G)$ for convenient vector spaces. Applications
to manifolds of mappings are given: The group of
$C^Q$-diffeomorphisms is a regular $C^Q$-Lie group but not better. \end{abstract}

\maketitle

Classes  of  Denjoy-Carleman  differentiable  functions  are in
general  situated  between  real  analytic functions and smooth
functions.  They  are  described  by  growth  conditions on the
derivatives.  Quasianalytic  classes  are  those where infinite
Taylor expansion is an injective mapping.
That  a  class  of  mappings  $\mathcal  S$ admits a convenient
setting  means  essentially  that  we  can  extend the class to
mappings   between   admissible   infinite  dimensional  spaces
$E,F,\dots$  so  that $\mathcal S(E,F)$ is again admissible and
we  have  $\mathcal  S(E\times  F,  G)$  canonically  $\mathcal
S$-diffeomorphic   to   $\mathcal  S(E,\mathcal  S(F,G))$  (the
exponential  law). Usually this comes hand in hand with (partly
nonlinear)   uniform   boundedness   theorems  which  are  easy
$\mathcal S$-detection principles.

For  the  $C^\infty$ convenient setting one can test smoothness
along smooth curves.
For  the real analytic ($C^\omega$) convenient setting we have:
A  mapping is $C^\omega$ if and only if it is $C^\infty$ and in
addition  $C^\omega$  along $C^\omega$-curves ($C^\omega$ along
just  affine  lines suffices). We shall use convenient calculus
of  $C^\infty$  and  $C^\omega$ mappings in this paper; see the
book  \cite{KM97}, or the three appendices in \cite{KMRc} for a
short overview.

In  \cite{KMRc}  we  succeeded  to  show that non-quasianalytic
log-convex  Denjoy-Carleman  classes  $C^M$  of moderate growth
(hence derivation closed) admit a convenient setting, where the
underlying admissible locally convex vector spaces are the same
as for smooth or for real analytic mappings. A mapping is $C^M$
if  and  only if it is $C^M$ along all $C^M$-curves. The method
of  proof  there relies on the existence of $C^M$ partitions of
unity.
In this paper we succeed to prove that quasianalytic log-convex
Denjoy-Carleman classes $C^Q$ of moderate growth which are also
$\mathcal  L$-intersections (see \thetag{\ref{nmb:1.6}}), admit
a convenient setting. The method consists of representing $C^Q$
as  the intersection $\bigcap \{C^L:L\in\mathcal L(Q)\}$ of all
larger  non-quasianalytic log-convex classes $C^L$; this is the
meaning   of:   $Q$   is   an   $\mathcal  L$-intersection.  In
\thetag{\ref{nmb:1.9}} we construct
countably   many   classes   $Q$   which   satisfy   all  these
requirements.
Taking  intersections  of derivation closed classes $C^L$ only,
or  only of classes $C^L$ of moderate growth, is not sufficient
for  yielding  the intended results. Thus we have to strengthen
many  results  from \cite{KMRc} before we are able to prove the
exponential  law. A mapping is $C^Q$ if and only if it is $C^L$
along each $C^L$-curve for each $L\in \mathcal L(Q)$.
It is an open problem (even in $\mathbb R^2$), whether a smooth
mapping which is $C^Q$ along each $C^Q$-curve (or affine line),
is indeed $C^Q$. As replacement we show that a mapping is $C^Q$
if  it  is  $C^Q$  along  each $C^Q$ mapping from a Banach ball
\thetag{\ref{nmb:5.2}}.
The  real  analytic  case  from  \cite{KrieglMichor90}  is  not
covered by this approach.

The  initial  motivation of both \cite{KMRc} and this paper was
the  desire  to  prove  the  following  result  which is due to
Rellich \cite{Rellich42V} in the real analytic case.
{\it  Let  $t\mapsto  A(t)$  for $t\in \mathbb R$ be a curve of
unbounded self-adjoint operators in a Hilbert space with common
domain  of  definition and with compact resolvent. If $t\mapsto
A(t)$ is of a certain quasianalytic
Denjoy-Carleman class $C^Q$, then
the   eigenvalues   and  the  eigenvectors  of  $A(t)$  may  be
parameterized $C^Q$ in $t$ also.
}
We  manage  to  prove this with the help of the results in this
paper  and in \cite{KMRc}. Due to length this will be explained
in another paper \cite{KMRp}.

Generally,  one  can  hope  that  the  space  $C^M(A,B)$ of all
Denjoy-Carleman   $C^M$-mappings   between  finite  dimensional
$C^M$-manifolds  (with  $A$  compact for simplicity) is again a
$C^M$-manifold,  that  composition is $C^M$, and that the group
$\on{Diff}^M(A)$ of all $C^M$-diffeomorphisms of $A$ is a
regular  infinite  dimensional  $C^M$-Lie group, for each class
$C^M$    which   admits   a   convenient   setting.   For   the
non-quasianalytic  classes  this was proved in \cite{KMRc}. For
quasianalytic classes this is proved in this paper.

\section{\label{nmb:1} Weight Sequences and function spaces}

\subsection{\label{nmb:1.1}Denjoy--Carleman $C^M$-functions in finite dimensions}

We mainly follow \cite{KMRc} and \cite{Thilliez08} (see also the references therein).
We use $\mathbb{N} = \mathbb{N}_{>0} \cup \{0\}$.
For each multi-index $\alpha=(\alpha_1,\ldots,\alpha_n) \in \mathbb{N}^n$, we write
$\alpha!=\alpha_1! \cdots \alpha_n!$, $|\alpha|= \alpha_1 +\cdots+ \alpha_n$, and 
$\partial^\alpha=\partial^{|\alpha|}/\partial x_1^{\alpha_1} \cdots \partial x_n^{\alpha_n}$. 

Let $M=(M_k)_{k \in \mathbb{N}}$ be a sequence of positive real numbers.
Let $U \subseteq \mathbb{R}^n$ be open.
We denote by $C^M(U)$ the set of all $f \in C^\infty(U)$ such that, for all
compact $K \subseteq U$,
there exist positive constants $C$ and $\rho$ such that
\begin{equation*}
|\partial^\alpha f(x)| \le C \, \rho^{|\alpha|} \, |\alpha|! \, M_{|\alpha|}
\quad\text{ for all }\alpha \in \mathbb{N}^n\text{ and }x \in K.
\end{equation*}
The set $C^M(U)$ is a \emph{Denjoy--Carleman class} of functions on $U$.
If $M_k=1$, for all $k$, then $C^M(U)$ coincides with the ring $C^\omega(U)$
of real analytic functions
on $U$. 
A sequence $M=(M_k)$ is \emph{log-convex} if $k\mapsto\log(M_k)$ is convex, i.e.,
\begin{equation*} M_k^2 \le M_{k-1} \, M_{k+1} \quad \text{ for all } k.
\end{equation*}

If $M=(M_k)$ is log-convex, then $k\mapsto (M_k/M_0)^{1/k}$ is increasing %
and
\begin{equation*}
M_l \, M_k\le M_0\,M_{l+k} \quad \text{ for all }l,k\in \mathbb{N}.
\tag{1}
\end{equation*}
Furthermore, we have that $k\mapsto k!M_k$ is log-convex (since
Euler's  $\Gamma$-function  is  so),  and  we  call this weaker
condition {\it weakly log-convex}.
If $M$ is weakly log-convex then $C^M(U,\mathbb R)$ is a ring, for all open
subsets $U \subseteq \mathbb{R}^n$. 
If $M$ is log-convex then (see the proof of \cite[2.9]{KMRc}) we have \begin{equation*}
M_1^j \, M_k\ge M_j\, M_{\alpha_1} \cdots M_{\alpha_j} \quad \text{ for all }
\alpha_i\in \mathbb{N}_{>0}\text{ with } \alpha_1+\dots+\alpha_j = k.
\tag{2}
\end{equation*}
This  implies  that the class of $C^M$-mappings is stable under
composition  (\cite{Roumieu62/63},  see  also \cite[4.7]{BM04};
this also follows from \thetag{\ref{nmb:1.4}}).
If $M$ is log-convex then the inverse function theorem
for $C^M$ holds (\cite{Komatsu79}; see also \cite[4.10]{BM04}),
and $C^M$ is closed under solving ODEs
(due to \cite{Komatsu80}).

Suppose that $M=(M_k)$ and $N=(N_k)$ satisfy $M_k \le C^k \, N_k$, for a constant $C$ and all $k$.
Then $C^M(U) \subseteq C^N(U)$. The converse
is true if $M$ is weakly log-convex:
There  exists  $f  \in C^M(\mathbb{R})$ such that $|f^{(k)}(0)|
\ge k! \, M_k$ for all $k$ (see \cite[Theorem 1]{Thilliez08}).

If  $M$  is  weakly  log-convex then $C^M$ is {\it stable under
derivations}  (alias  {\it  derivation  closed}) if and only if
\begin{equation*}\tag{3}
\sup_{k \in \mathbb{N}_{>0}} \Big(\frac{M_{k+1}}{M_k}\Big)^{\frac{1}{k}} < \infty.
\end{equation*}
A weakly log-convex sequence $M$ is called of \emph{moderate growth} if
\begin{equation*} \tag{4}
\sup_{j,k \in \mathbb{N}_{>0}} \Big(\frac{M_{j+k}}{M_j \, M_k}\Big)^{\frac{1}{j+k}} <
\infty.
\end{equation*}
Moderate growth implies derivation closed. 
\subsection*{Definition}

A  sequence  $M=(M_k)_{k=0,1,2,\dots}$  is called a {\it weight
sequence}  if  it  satisfies  $M_0=1\le M_1$ and is log-convex.
Consequently, it is increasing (i.e.\ $M_k\le M_{k+1}$).

A  {\it  DC-weight  sequence}  $M=(M_k)_{k=0,1,2,\dots}$  is  a
weight sequence which is also derivation closed
(DC stands for Denjoy-Carleman and also for derivation closed).
This was the notion investigated in \cite{KMRc}.

\begin{theorem}[Denjoy--Carleman \cite{Denjoy21}, \cite{Carleman26}]\label{nmb:1.2}
For a sequence $M$ of positive numbers the following statements are equivalent.
\begin{enumerate}
\item[(1)] $C^M$ is quasianalytic, i.e., for open connected 
$U \subseteq \mathbb{R}^n$ and each $a \in U$, the Taylor series homomorphism centered at $a$ 
from $C^M(U,\mathbb R)$ into the space of formal power series is injective.
\item[(2)] $\sum_{k=1}^\infty \frac1{m_k^{\flat(i)}} = \infty$
where $m_k^{\flat(i)}:= \inf\{(j!\,M_j)^{1/j}: j\ge k\}$ is the
increasing minorant of $(k!\,M_k)^{1/k}$.
\item[(3)] $\sum_{k=1}^\infty (\frac1{M_k^{\flat(lc)}})^{1/k} = \infty$
where $M^{\flat(lc)}_k$ is the log-convex minorant of $k!\,M_k$, given by 
$M^{\flat(lc)}_k := \inf\{(j!\,M_j)^{\frac{l-k}{l-j}}(l!\,M_l)^{\frac{k-j}{l-j}}: j\le k\le l, j<l\}$. 
\item[(4)] $\sum_{k=0}^\infty \frac{M^{\flat(lc)}_k}{M^{\flat(lc)}_{k+1}}=\infty$.
\end{enumerate}
\end{theorem}
For contemporary proofs see for instance \cite[1.3.8]{Hoermander83I} or \cite[19.11]{Rudin87}. 
\subsection{\label{nmb:1.3} Sequence spaces}
Let $M=(M_k)_{k\in\mathbb N}$ be a sequence of positive numbers and $\rho>0$.
We consider (where $\mathcal F$ stands for `formal power series')
\begin{align*}
\mathcal F^M_\rho &:= \Bigl\{(f_k)_{k\in\mathbb N}\in\mathbb R^\mathbb N:
\exists C>0\,\forall k\in\mathbb N:|f_k|\leq C\,\rho^k\;k!\;M_k\Bigr\}
\text{ and }
\mathcal F^M := \bigcup_{\rho>0}\mathcal F^M_\rho.
\end{align*}
Note  that, for $U\subseteq \mathbb R^n$ open, a function $f\in
C^\infty(U,\mathbb R)$ is in $C^M(U,\mathbb R)$ if and only if
for each compact $K \subset U$
\begin{equation*}
(\sup\{|\partial^\alpha f(x)|: x\in K, |\alpha|=k\})_{k\in \mathbb N} \in \mathcal F^M.
\end{equation*}

\begin{lemma*} We have
\begin{align*}
\mathcal F^{M^1}\subseteq \mathcal F^{M^2}
&\Leftrightarrow \exists\rho>0 \;\forall k: M^1_k\leq \rho^{k+1}\,M^2_k
\\&
\Leftrightarrow \exists C,\rho>0\; \forall k: M^1_k\leq C\,\rho^k\,M^2_k.
\end{align*}
\end{lemma*}
\begin{demo}{Proof}
($\Rightarrow$)
Let $f_k:=k!M^1_k$. Then
$f=(f_k)_{k\in\mathbb N}\in\mathcal F^{M^1}\subseteq \mathcal F^{M^2}$, so there
exists a $\rho>0$ such that $k!M^1_k\leq \rho^{k+1}k!M^2_k$ for all $k$.

($\Leftarrow$)
Let $f=(f_k)_{k\in\mathbb N}\in\mathcal F^{M^1}$, i.e.\
there exists a $\sigma>0$ with $|f_k|\leq\sigma^{k+1}k!M^1_k\leq (\rho\sigma)^{k+1}k!M^2_k$ for all $k$
and thus $f\in \mathcal F^{M^2}$.
\qed\end{demo}

\begin{lemma}\label{nmb:1.4}
Let $M$ and $L$ be sequences of positive numbers.
Then for the composition of formal power series we have
\[
\mathcal F^M\circ \mathcal F^L_{>0}\subseteq \mathcal F^{M\circ L}
\]
where 
$(M\circ L)_k := \max\{M_jL_{\alpha_1}\dots L_{\alpha_j}: 
\alpha_i\in \mathbb{N}_{>0}, \alpha_1+\dots+\alpha_j = k \}$ 
\end{lemma}

Here  $\mathcal  F^L_{>0} :=\{(g_k)_{k\in\mathbb N}\in \mathcal
F^L:  g_0=0\}$ is the space of formal power series in $\mathcal
F^L$ with vanishing constant term.

\demo{Proof} Let $f \in \mathcal F^M$ and $g \in \mathcal F^L$.
For $k>0$ we have (inspired by \cite{FaadiBruno1855})
\begin{align*}
\frac{(f\circ g)_k}{k!} &= \sum_{j=1}^k \frac{f_j}{j!}
\sum_{\substack{\alpha\in \mathbb{N}_{>0}^j\\ \alpha_1+\dots+\alpha_j =k}}
\frac{g_{\alpha_1}}{\alpha_1!}\dots
\frac{g_{\alpha_j}}{\alpha_j!}
\\
\frac{|(f\circ g)_k|}{k!(M\circ L)_k} &\le \sum_{j=1}^k \frac{|f_j|}{j!M_j}
\sum_{\substack{\alpha\in \mathbb{N}_{>0}^j\\ \alpha_1+\dots+\alpha_j =k}}
\frac{|g_{\alpha_1}|}{\alpha_1!L_{\alpha_1}}\dots
\frac{|g_{\alpha_j}|}{\alpha_j!L_{\alpha_j}}
\\&
\le \sum_{j=1}^k \rho_f^{j}C_f
\sum_{\substack{\alpha\in \mathbb{N}_{>0}^j\\ \alpha_1+\dots+\alpha_j =k}}
\rho_g^k C_g^j
\le \sum_{j=1}^k \rho_f^{j}C_f
\binom{k-1}{j-1}
\rho_g^k C_g^j
\\&
= \rho_g^k\rho_f C_f C_g \sum_{j=1}^k (\rho_f C_g)^{j-1}
\binom{k-1}{j-1}
= \rho_g^k\rho_f C_f C_g (1+\rho_f C_g)^{k-1}
\\&
= (\rho_g(1+\rho_f C_g))^k\frac{\rho_f C_f C_g}{1+\rho_f C_g}\qed
\end{align*}
\enddemo

\subsection{\label{nmb:1.5} Notation for quasianalytic weight sequences}
Let $M$ be a sequence of positive numbers. We may replace $M$ by $k\mapsto C\,\rho^k\,M_k$ with $C,\rho>0$
without changing $\mathcal F^M$. In particular, it is no loss of generality
to assume that $M_1>1$ (put $C\rho>1/M_1$) and $M_0=1$ (put $C:=1/M_0$).
If $M$ is log-convex then so is the modified sequence and if
in addition $\rho\geq M_0/M_1$ then the modified sequence is monotone increasing.
Furthermore  $M$  is  quasianalytic if and only if the modified
sequence is so, since $M_k^{\flat(lc)}$ is modified in the same
way.
We   tried  to  make  all  conditions  equivariant  under  this
modification.  Unfortunately,  the  next  construction does not
react nicely to this modification.

For a quasianalytic sequence $M=(M_k)$ let the sequence $\check M = (\check M_k)$
be defined by \begin{equation*}
\check M_k := M_k \prod_{j=1}^k \left(1-\frac1{(j!\,M_j)^{1/j}} \right)^k,\quad \check M_0=1.
\end{equation*}
We have $\check M_k\leq M_k$.
Note that if we put $m_k:=(k!M_k)^{1/k}$ (and $m_0:=1$) and 
$\check m_k:=(k!\check M_k)^{1/k}$ (where we assume $\check M_k \ge 0$)
then \[
\check m_k=m_k\prod_{j=1}^k \left(1-\frac1{m_j}\right)
\]
or, recursively, \[
\check m_{k+1}=\check m_k\frac{m_{k+1}-1}{m_k}\text{ and }\check m_0=1,\;\check m_1=m_1-1.
\]
And conversely, if all $\check M_k>0$ (this is the case if $M$ is increasing and $M_1>1$) then
\[
m_{k+1}=1+m_k\frac{\check m_{k+1}}{\check m_k}\text{ and }m_0=1,\;m_1=\check m_1+1
\]
i.e.\
\begin{equation*}\tag{1}
m_k=\check m_k\Bigl(1+\sum_{j=1}^k \frac1{\check m_j}\Bigr).
\end{equation*}

For sequences $M$ we define (recall from \thetag{\ref{nmb:1.1}}
that   $M$   is   called   weakly   log-convex   if   
$k\mapsto \log(k!\,M_k)$ is convex):
\begin{align*}
\mathcal L(M) &:= \{L\geq M: L\text{ non-quasianalytic, log-convex}\} \\
\mathcal L_w(M) &:= \{L\geq M: L\text{ non-quasianalytic, weakly log-convex}\}\supseteq \mathcal L(M) \\
\end{align*}

\begin{theorem}\label{nmb:1.6}
Let $Q=(Q_k)_{k=0,1,2,\dots}$ be a quasianalytic sequence of positive real numbers. Then we have:
\begin{enumerate}
\item[(1)] If %
the sequence $\check Q=(\check Q_k)$ is log-convex and positive then
\[
\mathcal F^Q = \bigcap_{L\in\mathcal L(Q)}\mathcal F^L.
\]
\item[(2)] 
If  $Q$ is weakly log-convex, then for each $L^1,L^2\in\mathcal
L_w(Q)$  there  exists  an  $L\in  \mathcal L_w(Q)$ with $L \le L^1,L^2$.
\item[(3)]
If $Q$ is weakly log-convex of moderate growth, then
for each $L\in\mathcal L_w(Q)$ there exists an $L'\in\mathcal L_w(Q)$ such that 
$L'_{j+k} \le C^{j+k} L_j L_k$ for some positive constant $C$ and all $j,k \in \mathbb N$.
\end{enumerate}
\end{theorem}
We could not obtain \thetag{2} for log-convex instead of weakly
log-convex,  in  particular  for  $\mathcal  L(Q)$  instead  of
$\mathcal L_w(Q)$.

\subsection*{Definition}   A   quasianalytic  sequence  $Q$  of
positive    real    numbers    is    called    {\it   $\mathcal
L$-intersectable}  or  an  {\it  $\mathcal  L$-intersection} if
$\mathcal F^Q = \bigcap_{L\in\mathcal L(Q)}\mathcal F^L$ holds.
Note that we may replace any non-quasianalytic weight sequence $L$ for which
$k\mapsto (\frac{Q_k}{L_k})^{1/k}$ is bounded, by an $\tilde L\in \mathcal L(Q)$
with $\mathcal F^{\tilde L}=\mathcal F^L$:
Choose %
$\rho\geq 1/L_1$ (see \thetag{\ref{nmb:1.5}}) and $\rho\geq \sup\{(\frac{Q_k}{L_k})^{1/k}:k\in\mathbb N\}$
then $\tilde L_k:=\rho^k L_k\geq Q_k$. 
\begin{demo}{Proof}
\thetag{1} The proof is partly adapted from \cite{Boman63}.

Let $q_k = (k!\,Q_k)^{1/k}$ and $q_0=1$, similarly 
$\check q_k = (k!\,\check Q_k)^{1/k}$, $l_k=(k!L_k)^{1/k}$, etc.
Then $\check q$ is increasing since $\check Q_0=1$, and $\check Q$
and the Gamma function are log-convex.

Clearly $\mathcal F^Q \subseteq \bigcap_{L\in\mathcal L(Q)}\mathcal F^L$. 
To show the converse inclusion, let $f\notin \mathcal F^Q$ and $g_k:=|f_k|^{1/k}$. %
Then
\[
\varlimsup \frac{g_k}{q_k} = \infty.
\]
Choose $a_j, b_j>0$ with $a_j\nearrow \infty$, $b_j\searrow 0$, and $\sum\frac1{a_jb_j}<\infty$.
There exist strictly increasing $k_j$ such that $\frac{g_{k_j}}{q_{k_j}}\ge a_j$. 
Since $\frac{q_k}{\check q_k}$ is increasing by \thetag{\ref{nmb:1.5}.1} we get
$b_j\,\frac{g_{k_j}}{\check q_{k_j}}= b_j\,\frac{g_{k_j}}{q_{k_j}}\,\frac{q_{k_j}}{\check q_{k_j}}
\ge a_j b_j\,\frac{q_{k_1}}{\check q_{k_1}}\to \infty$.
Passing to a subsequence we may assume that $k_0>0$ and 
$1 < \beta_j:=b_j\,\tfrac{g_{k_j}}{\check q_{k_j}}\nearrow \infty$. 
Passing to a subsequence again we may also get 
\begin{equation*} \tag{4}
\beta_{j+1} \ge \left(\beta_j\right)^{k_j}.
\end{equation*}

Define a piecewise affine function $\phi$ by \begin{align*}
\phi(k) := \begin{cases}
0 & \text{ if } k=0,
\\
k_j \log\beta_j & \text{ if } k=k_j,
\\
c_j + d_j k & \text{ for the minimal $j$ with } k \le k_j,
\end{cases}
\end{align*}
where $c_j$ and $d_j$ are chosen such that $\phi$ is well defined and 
$\phi(k_{j-1})=c_j+d_j k_{j-1}$, i.e., for $j\ge 1$,
\begin{align*} \tag{5}
c_j+d_j k_j &= k_j \log\beta_j,\\
c_j+d_j k_{j-1} &= k_{j-1} \log \beta_{j-1},
\quad \text{ and } \\
c_0 &=0,\\
d_0 &= \log\beta_0.
\end{align*}
This implies first that $c_j\le 0$ and then
\begin{align*} \tag{6}
\log\beta_j \le d_j&=\frac{k_j\log \beta_j-k_{j-1}\log\beta_{j-1}}{k_j-k_{j-1}} 
\le \frac{k_j}{k_j-k_{j-1}} \log\beta_j
\\&
\overset{(4)}\leq\frac{\log\beta_{j+1}}{k_j-k_{j-1}}\leq\log\beta_{j+1}.
\end{align*}
Thus $j\mapsto d_j$ is increasing.
It follows that $\phi$ is convex. The fact that all $c_j \le 0$ implies that $\phi(k)/k$ is increasing.

Now let
\[
L_k := e^{\phi(k)} \cdot \check Q_k.
\]
Then   $L=(L_k)$   is   log-convex  and  satisfies  $L_0=1$  by
construction   and   $f\notin  \mathcal  F^L$,  since  we  have
$\frac{l_{k_j}}{g_{k_j}}=\frac{\check q_{k_j}\beta_j}{g_{k_j}}=
b_j \to 0$ and so $\varlimsup \frac{g_k}{l_k}= \infty$.
Let  us  check that $L$ is not quasianalytic. By \thetag{6} and
since $(\check q_k)$ is increasing, we have, for 
$k_{j-1} \le k < k_j$,
\begin{align*}
\frac{L_k}{(k+1)L_{k+1}} &= \frac{e^{\phi(k)-\phi(k+1)}\,\check Q_k}{(k+1)\,\check Q_{k+1}}
= \frac{e^{\phi(k)-\phi(k+1)} \,\check q_k^k}{\check q_{k+1}^{k+1}}
=e^{-d_j} \frac{\check q_k^k}{\check q_{k+1}^{k+1}} \\&
\le \frac{1}{\beta_j\,\check q_k}
=\frac{\check q_{k_j}}{b_j g_{k_j}} \frac{1}{\check q_k}.
\end{align*}
Thus, by \thetag{\ref{nmb:1.5}.1},
\[
\sum_{k=k_{j-1}}^{k_j-1} \frac{L_k}{(k+1)L_{k+1}} 
\le \frac{\check q_{k_j}}{b_j g_{k_j}} \sum_{k=k_{j-1}}^{k_j-1} \frac{1}{\check q_k}
\le \frac{q_{k_j}}{b_j g_{k_j}} \le \frac{1}{a_j b_j}, 
\]
which shows that $L$ is not quasianalytic and $C_1:=\sum_{k=1}^\infty \tfrac1{l_k} < \infty$
by \thetag{\ref{nmb:1.2}}.

Next we claim that $\mathcal F^Q \subseteq \mathcal F^L$. 
Since $\tfrac{l_k}{\check q_k} = \frac{(k!L_k)^{1/k}}{(k!\check Q_k)^{1/k}} =e^{\phi(k)/k}$ 
is increasing, we have
\[
\infty>\frac{\check q_1}{l_1}+C_1 > \frac{\check q_1}{l_1} + \sum_{j=1}^k \frac1{l_j} 
= \frac{\check q_1}{l_1} + \sum_{j=1}^k \frac{\check q_j}{l_j} \frac1{\check q_j} 
\ge \frac{\check q_k}{l_k} \Big(1+\sum_{j=1}^k \frac1{\check q_j}\Big)
= \frac{q_k}{l_k},
\]
which proves $\mathcal F^Q \subseteq \mathcal F^L$.
Finally we may replace $L$ by some $L\in\mathcal L(Q)$ without changing $\mathcal F^L$
by the remark before the proof.
Thus \thetag{1} is proved. 
\thetag{2} Assume without loss that $L^1_0=L^2_0=1$.
Let %
$k!L_k$ be the log-convex minorant of $k!\bar L_k$ where $\bar L_k:=\min\{L^1_k,L^2_k\}$. 
Since $L^1,L^2\geq \bar L \ge Q$ and $k!Q_k$ is log-convex we have $L^1,L^2\geq L\geq Q$. 
Since $L^1, L^2$ are not quasianalytic and are weakly log-convex (hence
$k\mapsto (k!L_k^j)^{1/k}$ is increasing), we get that $k\mapsto (k!\bar L_k)^{1/k}$
is increasing and
\begin{align*}
\sum_k \frac{1}{(k!\,\bar L_k)^{1/k}} \le \sum_k \frac{1}{(k!\,L^1_k)^{1/k}} + \sum_k \frac{1}{(k!\,L^2_k)^{1/k}}
< \infty.
\end{align*}
By (\ref{nmb:1.2}, 2$\Rightarrow$1) we get that ${\bar L}$
is not quasianalytic. 
By (\ref{nmb:1.2}, 1$\Rightarrow$3) we get $\sum_k \frac1{(k!L_k)^{1/k}}<\infty$ since $\bar L^{\flat(lc)}=L$, i.e.\
$L$ is not quasianalytic.

\thetag{3}
Let $\tilde Q_k := k! Q_k$, $\tilde L_k := k! L_k$, and so on.
Since $Q$ is of moderate growth we have 
\[
C_{\tilde Q}:= \sup_{k,j}\left(\frac{\tilde Q_{k+j}}{\tilde Q_k \tilde Q_j}\right)^{1/(k+j)} 
\le 2 \sup_{k,j}\left(\frac{Q_{k+j}}{Q_k Q_j}\right)^{1/(k+j)}<\infty.
\]
Let $L\in \mathcal L_w(Q)$; without loss we assume that $L_0=1$. 
We put \begin{align*}
\tilde L'_k :&= C_{\tilde Q}^k\min\{\tilde L_j \tilde L_{k-j}: j=0,\dots ,k\}
= C_{\tilde Q}^k\min\{\tilde L_j \tilde L_{k-j}: 0\le j\le k/2\}.
\end{align*}
Then
\[
\sup_{k,j}\left(\frac{L'_{k+j}}{L_k L_j}\right)^{1/(k+j)} 
\le \sup_{k,j}\left(\frac{\tilde L'_{k+j}}{\tilde L_k \tilde L_j}\right)^{1/(k+j)} 
\le C_{\tilde Q} < \infty.
\]
Since $\tilde L$ is log-convex we have
$\tilde L_k^2\le \tilde L_j \tilde L_{2k-j}$ and 
$\tilde L_k \tilde L_{k+1}\le \tilde L_j \tilde L_{2k+1-j}$ for $j=0,\dots, k$;
therefore
$\tilde L'_{2k} = C_{\tilde Q}^{2k}\tilde L_k^2$ and 
$\tilde L'_{2k+1} = C_{\tilde Q}^{2k+1}\tilde L_k \tilde L_{k+1}$. 
It is easy to check that $\tilde L'$ %
is log-convex. 
To see that $L'$ is not quasianalytic we will use
that $(\tilde L'_k)^{1/k}$ is increasing since $\tilde L'$ is log-convex. 
So it suffices to compute the sum of the even indices only.
\begin{align*}
\sum_{k}\frac1{\tilde L'_{2k}{}^{1/(2k)}}
= \frac1{C_{\tilde Q}}\sum_k\frac1{\tilde L_k{}^{1/k}} < \infty.
\end{align*}
It remains to show that $L' \ge Q$.
Since $L\in \mathcal L_w(Q)$ we have $Q\le L$ and for $j=\lfloor k/2\rfloor$,
\begin{align*}
\frac{Q_k}{L'_k} = \frac{\tilde Q_k}{\tilde L'_k} &= \frac{\tilde Q_k}{C_{\tilde Q}^k\, \tilde L_j \tilde L_{k-j}}
\le \frac{\tilde Q_k}{\frac{\tilde Q_k}{\tilde Q_j \tilde Q_{k-j}} \tilde L_j \tilde L_{k-j}}
\le \frac{\tilde Q_j}{\tilde L_j}\frac{\tilde Q_{k-j}}{\tilde L_{k-j}} \le1.
\qed\end{align*}
\end{demo}

\begin{corollary}\label{nmb:1.7}
Let $Q$ be a quasianalytic weight sequence. 
Then \[
\mathcal F^Q = \bigcap _{L\in \mathcal L_w(Q)} \mathcal F^L.
\]
\end{corollary}

\begin{demo}{Proof}
Without loss we may assume that the sequence $\check q_k$ is increasing.
Namely, by definition this is the case if and only if $q_k\le q_{k+1}-1$. 
Since $Q_0=1$ and $(Q_k)$ is log-convex, $Q_k^{1/k}$ is increasing and thus
$q_{k+1}-q_k \ge Q_k^{\frac1k}((k+1)!^{\frac1{k+1}} - k!^{\frac1k}) \ge Q_1\,\frac 1e \ge \frac 1e$. 
If we set $\tilde Q_k := e^k Q_k$, then $\tilde Q=(\tilde Q_k)$ 
is a quasianalytic weight sequence with 
$\tilde Q_1>1$, $\mathcal F^{\tilde Q} = \mathcal F^Q$, and 
$\check {\tilde q}_k$ is increasing.

Now a little adaptation of the proof of \thetag{\ref{nmb:1.6}.1} shows the corollary:
Define here
\[
l_k := \beta_j \check q_k\quad \text{ for the minimal }j\text{ with } k\le k_j.
\]
Then $\frac{l_{k_j}}{g_{k_j}}= \frac{\beta_j \check q_{k_j}}{g_{k_j}}= b_j \to 0$ 
and so $\varlimsup \frac{g_k}{l_k}= \infty$.
We have
\[
\sum_{k=k_{j-1}+1}^{k_j} \frac1{l_k} = \sum_{k=k_{j-1}+1}^{k_j} \frac{1}{\beta_j\check q_k}
= \frac{\check q_{k_j}}{b_{j}g_{k_j}}\sum_{k=k_{j-1}+1}^{k_j} \frac1{\check q_k}
\le \frac{q_{k_j}}{b_{j}g_{k_j}}
\le \frac1{a_jb_j}
\]
and thus $\sum_{k=1}^\infty \frac1{l_k}< \infty$.
As $l_k$ is increasing, the Denjoy--Carleman theorem \thetag{\ref{nmb:1.2}} 
implies that $L_k=\frac{l_k^k}{k!}$ is non-quasianalytic.
Since $\frac{l_k}{\check q_k} = \beta_j$ is increasing,
we find (as in the proof of \thetag{\ref{nmb:1.6}.1}) that $C:=\max\{L_0/L_1,\sup_k \frac{q_k}{l_k}\} <\infty$. 
Replacing $L_k$ by $C^k L_k$ we may assume that $Q\le L$.
Let the sequence $k!\underline L_k$ be the log-convex minorant of $k! L_k$.
Since $Q_k$ is (weakly) log-convex, we have $Q \le \underline L$.
By \thetag{\ref{nmb:1.2}} and the fact that $L$ is non-quasianalytic, $\underline L$ 
is non-quasianalytic as well.
Thus $\underline L \in \mathcal L_w(Q)$ and still $f\notin \mathcal F^{\underline L}$.
\qed\end{demo}

Corollary \thetag{\ref{nmb:1.7}} implies that for the sequence $\omega=(1)_k$ 
describing real analytic functions we have 
$\mathcal F^\omega = \bigcap_{L\in \mathcal L_w(\omega)} \mathcal F^L$. 
Note that $\mathcal L_w(\omega)$ consists of all weakly log-convex non-quasianalytic $L\ge 1$. 
This is slightly stronger than a result by T.\ Bang, who shows that 
$\mathcal F^\omega = \bigcap\mathcal F^L$ where $L$ runs through all 
non-quasianalytic sequences with $l_k=(k!\,L_k)^{1/k}$ increasing, see \cite{Bang46}, \cite{Boman63}.

This result becomes wrong if we replace weakly log-convex by log-convex:

\subsection{\label{nmb:1.8}The intersection of all $\mathcal F^L$, where $L$ is any 
non-quasianalytic weight sequence} Put \[
Q_k := \frac{(k \log(k+e))^k}{k!}, \quad Q_0 :=1.
\] Then $Q=(Q_k)$ is a quasianalytic weight sequence of moderate growth with $Q_1>1$.
We claim that $Q$ is $\mathcal L$-intersectable, i.e., 
$\mathcal F^Q = \bigcap_{L\in\mathcal L(Q)}\mathcal F^L$. 
We could check that $\check Q$ is log-convex. 
This can be done, but is quite cumbersome.
A simpler argument is the following. 
We consider $\check q'_k:=k$, $\check q'_0:=1$. 
Then $\check Q_k' =k^k/k!$ is log-convex. 
Since $C_1 \log k \le \sum_{j=1}^k \frac{1}{j} \le C_2 \log k$, we have by \thetag{\ref{nmb:1.5}.1}
\begin{equation*} C_3 k \log(k+e) \le q'_k \le C_4 k \log(k+e)
\end{equation*}
for suitable constants $C_i$. 
Hence $\mathcal F^Q=\mathcal F^{Q'}$.
By theorem \thetag{\ref{nmb:1.6}.1} we have
\[
\mathcal F^Q=\mathcal F^{Q'} = \bigcap_{L \in \mathcal L(Q')} \mathcal F^L 
= \bigcap_{L \in \mathcal L(Q)} \mathcal F^L
\]
since $\mathcal L(Q)$ and $\mathcal L(Q')$ contain only sequences 
which are "equivalent mod $(\rho^k)$".
The claim is proved. 
Let $L$ be any non-quasianalytic weight sequence. Consider
\[
\alpha_k:= \frac{(k!\, L_k)^{\frac{1}{k}}}{k} = \frac{l_k}{k}.
\]
Since $L$ is log-convex and $L_0=1$, we find that $L_k^{1/k}$ is increasing. 
Thus, for $s \le k$ we find \[
\frac{\alpha_s}{\alpha_k} =\frac{k}{s} \cdot \frac{s!^{1/s}}{k!^{1/k}} \cdot \frac{L_s^{1/s}}{L_k^{1/k}} \le 2e
\]
(using Stirling's formula for instance).
Since $L$ is not quasianalytic, we have $\sum_{k=1}^\infty \frac{1}{k \alpha_k} < \infty$. 
But
\[
\sum_{\sqrt{k} \le s \le k} \frac{1}{s \alpha_s} 
\ge \frac{1}{2e} \cdot \frac{1}{\alpha_k} \sum_{\sqrt{k} 
\le s \le k} \frac{1}{s} \sim \frac{1}{2e} \cdot \frac{1}{\alpha_k} \cdot \frac{\log k}{2}.
\]
The sum on the left tends to $0$ as $k \to \infty$. 
So $\frac{\log k}{\alpha_k} = \frac{k \log k}{l_k}$
is bounded.
Thus $\mathcal F^{Q} \subseteq \mathcal F^L$.

So we have proved the following theorem (which is intimately related to \cite[Thm.\ C]{Rudin62}).

\begin{theorem*}
Put $Q_k=(k \log (k+e))^k/k!$, $Q_0=1$. 
Then $Q$ is $\mathcal L$-intersectable. 
In fact, \[
\mathcal F^Q = \bigcap\{\mathcal F^L : L \text{ non-quasianalytic weight sequence}\}.\qed
\]
\end{theorem*}

\subsection*{Remark}
Log-convexity of $\check Q$ is only sufficient for $Q$ being an 
$\mathcal L$-intersection, see \thetag{\ref{nmb:1.6}.1}:
Using Stirling's formula we see that $\mathcal F^{Q}=\mathcal F^{Q''}$ for 
$Q_k=(k \log (k+e))^k/k!$ and $Q''_k = (\log(k+e))^k$. 
Also $\mathcal L(Q)$ and $\mathcal L(Q'')$ contain only sequences which are 
``equivalent mod $(\rho^k)$'' and
\thetag{\ref{nmb:1.6}.1} holds for $Q$, thus also for $Q''$. 
But $\check Q''$ is not log-convex. 
\subsection{\label{nmb:1.9} A class of examples}

Let $\log^n$ denote the $n$-fold composition of $\log$ defined recursively by
\begin{align*}
\log^1 &:= \log,\\
\log^n &:= \log \circ \log^{n-1}, \quad (n\ge2).
\end{align*}
For $0< \delta \le 1$, $n \in \mathbb{N}_{>0}$, we recursively define sequences 
$q^{\delta,n}=(q^{\delta,n}_k)_{k \ge \kappa_n}$ by
\begin{align*}
q^{1,1}_k &:= k \log k,\\
q^{\delta,n}_k &:= q^{1,n-1}_k \cdot (\log^n(k))^{\delta}, \quad (n\ge2),
\end{align*}
where $\kappa_n$ is the smallest integer greater than $e \uparrow \uparrow n$, i.e.,
\[
\kappa_n := \lceil e \uparrow \uparrow n \rceil, \quad
e \uparrow \uparrow n := \underbrace{e^{e^{{\cdot}^{{\cdot}^ e}}}}_{n \text{ times}}. 
\]
Let $Q^{\delta,n}:=(Q^{\delta,n}_k)_{k \in \mathbb{N}}$ with \begin{align*}
Q^{\delta,n}_0 &:=1, \\
Q^{\delta,n}_k &:= \frac{1}{(k-1+\kappa_n)!}(q^{\delta,n}_{k-1+\kappa_n})^{k-1+\kappa_n}, \quad (k\ge1), 
\end{align*}
and consider
\[
\mathcal Q :=\{Q^{1,1}\}\cup \{Q^{\delta,n} : 0 < \delta \le 1, n \in \mathbb{N}_{>1}\}. 
\]
It is easy to check inductively that each $Q \in \mathcal Q$ 
is a quasianalytic weight sequence of moderate growth with $Q_1>1$.
Namely, $(\log^n(k))^{\delta k}$ is increasing, log-convex, and has moderate growth. 
Quasianalyticity follows from Cauchy's condensation criterion or the integral test.
By construction, $\mathcal Q\ni Q\mapsto \mathcal F^Q$ is injective.

Let us consider
\[
\hat q^{1,n}_k := q^{1,n-1}_k \left(1+\sum_{j=\kappa_n}^k \frac{1}{q^{1,n-1}_j}\right). 
\]
Since $\frac{d}{dx} \log^n (x) = \frac{1}{x \log (x) \cdots \log^{n-1} (x)}$, we have 
(by comparison with the corresponding integral)
\[
C_1 \log^n(k) \le \sum_{j=\kappa_n}^k \frac{1}{q^{1,n-1}_j} \le C_2 \log^n(k)
\]
and thus
\begin{equation*} \tag{1}
C_3 q^{1,n}_k \le \hat q^{1,n}_k \le C_4 q^{1,n}_k
\end{equation*}
for suitable constants $C_i$. 
Hence $\mathcal F^{Q^{1,n}} = \mathcal F^{\hat Q^{1,n}}$.
Since $Q^{1,n-1}$ is log-convex, theorem \thetag{\ref{nmb:1.6}.1} implies
\[
\mathcal F^{Q^{1,n}} = \mathcal F^{\hat Q^{1,n}} 
= \bigcap_{L \in \mathcal L(\hat Q^{1,n})} \mathcal F^L 
= \bigcap_{L \in \mathcal L(Q^{1,n})} \mathcal F^L
\]
since $\mathcal L(\hat Q^{1,n})$ and $\mathcal L(Q^{1,n})$ 
contain only sequences which are "equivalent mod $(\rho^k)$".

Hence we have proved (the case $n=1$ follows from \thetag{\ref{nmb:1.8}}):

\begin{theorem*}
Each $Q^{1,n}$ ($n \in \mathbb{N}_{>0}$) is a quasianalytic weight sequence of moderate growth 
which is an $\mathcal L$-intersection, i.e.,
\[
\mathcal F^{Q^{1,n}} = \bigcap_{L \in \mathcal L(Q^{1,n})} \mathcal F^L. \qed
\]
\end{theorem*}

\subsection*{Conjecture}
{\it This is true for each $Q \in \mathcal Q$.}

\subsection*{Remark}
Let $\check Q$ be any quasianalytic log-convex sequence of positive numbers. 
Then the corresponding sequence $Q$
(determined by \thetag{\ref{nmb:1.5}.1}) is quasianalytic and $\mathcal L$-intersectable.
However, the mapping $\check Q \mapsto \mathcal F^{Q}$ is not injective.
For instance, the image of $(C \rho^k \check Q_k)_k$ is the same for all positive $C$ and $\rho$
(which follows from \thetag{\ref{nmb:1.5}.1}).
Here is a more striking example:

Let $Q^{\delta,n} \in \mathcal Q$ and let $P^{\delta,n}=(P_k^{\delta,n})_k$ be defined by
\[
P_k^{\delta,n} := \frac{1}{(k-1+\kappa_n)!}(p^{\delta,n}_{k-1+\kappa_n})^{k-1+\kappa_n}, \quad P^{\delta,n}_0:=1,
\]
where
\begin{align*}
p^{\delta,n}_k &:= q^{\delta,n}_k \left(1+\sum_{j=\kappa_n}^k \frac{1}{q^{\delta,n}_j}\right), \quad 
\text{ for } 0<\delta <1,\\
p^{1,n}_k &:= \hat q^{1,n+1}_k = q^{1,n}_k \left(1+\sum_{j=\kappa_{n+1}}^k \frac{1}{q^{1,n}_j}\right).
\end{align*}
{\it We claim that \; $\mathcal F^{P^{1,n-1}} = \mathcal F^{P^{\delta,n}} = \mathcal F^{P^{\epsilon,n}}$ 
for all $0 < \delta, \epsilon < 1$.}
For:
Since
\[
\frac{d}{dx} \frac{(\log^n(x))^{1-\delta}}{1-\delta}
= \frac{1}{x \log (x) \cdots \log^{n-1} (x) (\log^n(x))^\delta},
\]
we have
\[
C_1 \frac{(\log^n(k))^{1-\delta}}{1-\delta} 
\le \sum_{j=\kappa_n}^k \frac{1}{q^{\delta,n}_j} 
\le C_2 \frac{(\log^n(k))^{1-\delta}}{1-\delta},
\]
and thus
\[
\frac{p^{\delta,n}_k}{p^{\epsilon,n}_k} 
= \frac{(\log^n(k))^\delta}{(\log^n(k))^\epsilon}
\frac{\Big(1+\sum_{j=\kappa_n}^k \frac{1}{q^{\delta,n}_j}\Big)}{\Big(1+\sum_{j=\kappa_n}^k 
\frac{1}{q^{\epsilon,n}_j}\Big)} \le
C_3 \frac{(\log^n(k))^\delta}{(\log^n(k))^\epsilon} \frac{(\log^n(k))^{1-\delta}}{(\log^n(k))^{1-\epsilon}} = C_3
\]
and similarly
\[
\frac{p^{\delta,n}_k}{p^{\epsilon,n}_k} \ge C_4
\]
for suitable constants $C_i$.
By lemma \thetag{\ref{nmb:1.3}} we have $\mathcal F^{P^{\delta,n}} = \mathcal F^{P^{\epsilon,n}}$ for all 
$0 < \delta, \epsilon < 1$.
The same reasoning with $\delta=0$ proves that 
$\mathcal F^{P^{1,n-1}} = \mathcal F^{P^{\epsilon,n}}$.

\subsection{\label{nmb:1.10} Definition of function spaces}
Let $M=(M_k)_{k\in\mathbb N}$ be a sequence of positive numbers,
$E$ and $F$ be Banach spaces, $U\subseteq E$ open, $K\subseteq U$
compact, and $\rho>0$.
We consider the non-Hausdorff Banach space
\begin{align*}
C^M_{K,\rho}(U,F) &:= \Bigl\{f\in C^\infty(U,F):(\sup_{x\in K}\|f^{(k)}(x)\|_{L^k(E,F)})_k\in\mathcal F^M_\rho\Bigr\} \\
&=\Bigl\{f\in C^\infty(U,F): \|f\|_{K,\rho} <\infty\Bigr\}, \quad \text{where}\\
\|f\|_{K,\rho} &:=\sup\Bigl\{\frac{\|f^{(k)}(x)\|_{L^k(E,F)}}{ k!\,M_k\,\rho^k}:x\in K,k\in\mathbb N\Bigr\},
\intertext{the inductive limit}
C^M_K(U,F)&:=\varinjlim_{\rho>0} C^M_{K,\rho}(U,F),\\
\intertext{and the projective limit} C_b^M(U,F)&:=\varprojlim_{K\subseteq U} C^M_K(U,F), 
\text{ where $K$ runs through all compact subsets of $U$}.
\end{align*}
Here $f^{(k)}(x)$ denotes
the $k^{\text{th}}$-order Fr\'echet derivative of $f$ at $x$.

Note that instead of $\|f^{(k)}(x)\|_{L^k(E,F)}$ we could equivalently use
$\sup\{\|d^k_vf(x)\|_F:\|v\|_E\leq 1\}$ by \cite[7.13.1]{KM97}. 
For $E=\mathbb R^n$ and $F=\mathbb R$ this is the same space as in \thetag{\ref{nmb:1.1}}.

For convenient vector spaces $E$ and $F$, and $c^\infty$-open $U\subseteq E$ we define:
\begin{align*}
C_b^M(U,F) &:= \Bigl\{f\in C^\infty(U,F):
\forall B\; \forall\text{ compact } K\subseteq U\cap E_B\;\exists \rho>0:\\
&\quad\qquad\big\{\frac{f^{(k)}(x)(v_1,\dots,v_k)}{k!\,\rho^k\, M_k}:k\in \mathbb N,x\in K,\|v_i\|_B\leq 1\big\}
\text{ is bounded in $F$}\Bigr\}
\\&
= \Bigl\{f\in C^\infty(U,F):
\forall B\; \forall\text{ compact } K\subseteq U\cap E_B\;\exists \rho>0:\\
&\quad\qquad\big\{\frac{d^k_vf(x)}{k!\,\rho^k\, M_k}:k\in \mathbb N,x\in K,\|v\|_B\leq 1\big\}\text{ is bounded in $F$}\Bigr\}.
\end{align*}
Here $B$ runs through all closed absolutely convex bounded subsets and $E_B$ is the vector space generated by $B$ with the Minkowski functional $\|v\|_B= \inf\{\lambda\ge 0: v\in \lambda B\}$ as complete norm.

Now we define the spaces of main interest in this paper: First we put
\[
C^M(\mathbb R,U) := \{c:\mathbb R\to U\;:\;\ell\circ c\in C^M(\mathbb R,\mathbb R)\;\forall \ell\in E^*\}. 
\]
In general, for $L$ log-convex non-quasianalytic we put
\begin{align*}
C^L(U,F) &:= \{f:f\circ c\in C^L(\mathbb R,F)\;\forall c\in C^L(\mathbb R,U)\}\\
&= \{f:\ell\circ f\circ c\in C^L(\mathbb R,\mathbb R)\;\forall c\in C^L(\mathbb R,U),\forall\ell\in F^*\}
\end{align*}
supplied with the initial locally convex structure induced by all linear mappings 
$C^L(c,\ell): f\mapsto \ell\circ f\circ c \in C^L(\mathbb R, \mathbb R)$, 
which is a convenient vector space as $c^\infty$-closed subspace in the product.
Note that in particular the family $\ell_*:C^L(U,F)\to C^L(U,\mathbb R)$ with $\ell\in F^*$ 
is initial, whereas this is not the case for $C^L$ replaced by $C_b^L$ as example \thetag{\ref{nmb:1.11}} for
$\{\on{inj}_k\circ g^\vee(k):k\in\mathbb N\}\subseteq C^L(\mathbb R,\mathbb R^\mathbb N)$ shows, where
$\on{inj}_k$ denotes the inclusion of the $k$-th factor in $\mathbb R^\mathbb N$.
\newline
For $Q$ a quasianalytic $\mathcal L$-intersection
we define the space
\begin{align*}
C^Q(U,F) &:= \bigcap_{L\in\mathcal L(Q)} C^L(U,F) %
\end{align*}
supplied with the initial locally convex structure.
By theorem \thetag{\ref{nmb:1.6}.1} this definition coincides with the classical notion of 
$C^Q$ if $E$ and $F$ are finite dimensional.

\begin{lemma*}
For $Q$ a quasianalytic $\mathcal L$-intersection, %
the composite of $C^Q$-mappings is again $C^Q$,
and bounded linear mappings are $C^Q$.
\end{lemma*}

\begin{demo}{Proof}
This is true for $C^L$ (see \cite[3.1 and 3.11.1]{KMRc}) for every 
$L\in\mathcal L(Q)$ since each such $L$ is log-convex.
\qed\end{demo}

\subsection{\label{nmb:1.11} Example}
By \cite[Theorem 1]{Thilliez08}, for each weakly log-convex sequence $M$ there exists 
$f \in C^M(\mathbb{R},\mathbb{R})$ such that $|f^{(k)}(0)| \ge k! \, M_k$ for all $k \in \mathbb{N}$. 
{\it Then $g:\mathbb R^2 \to \mathbb R$ given by $g(s,t) = f(st)$ is $C^M$, whereas there is no reasonable
topology on $C^M(\mathbb R,\mathbb R)$ such that the associated mapping 
$g^\vee:\mathbb R\to C^M(\mathbb R,\mathbb R)$ is $C_b^M$.}
For a topology on $C^M(\mathbb R,\mathbb R)$ to be
reasonable we require only that all evaluations $\on{ev}_t: C^M(\mathbb R,\mathbb R) \to \mathbb R$ are
bounded linear functionals.

\begin{demo}{Proof} The mapping $g$ is obviously $C^M$. 
If $g^\vee$ were $C_b^M$, for $s=0$ there existed $\rho$ such that \[
\Big\{\frac{(g^\vee)^{(k)}(0)}{k!\,\rho^k\,M_k}:k\in\mathbb N\Big\}
\]
was bounded in $C^M(\mathbb R,\mathbb R)$. 
We apply the bounded linear functional $\on{ev}_{t}$ for $t=2\rho$ and then get \[
\frac{(g^\vee)^{(k)}(0)(2\rho)}{k!\,\rho^k\,M_k}=\frac{(2\rho)^kf^{(k)}(0)}{k!\,\rho^k\,M_k} \ge 2^k,
\]
a contradiction.
\qed\end{demo}

This example shows that for $C_b^M$ one cannot expect cartesian closedness. 
Using cartesian closedness \thetag{\ref{nmb:3.3}} and \thetag{\ref{nmb:2.3}} this also shows 
(for $F=C^M(\mathbb R,\mathbb R)$ and $U=\mathbb R=E$) that
\[
C_b^M(U,F)\supsetneq\bigcap_{B,V} C_b^M(U\cap E_B,F_V)
\]
where $F_V$ is the completion of $F/p_V^{-1}(0)$ with respect to the seminorm 
$p_V$ induced by the absolutely convex closed 0-neighbourhood $V$.

If we compose $g^\vee$ with the restriction map $(\on{incl}_\mathbb N)^*:C^M(\mathbb R,\mathbb R)\to\mathbb R^\mathbb N:=\prod_{t\in\mathbb N}\mathbb R$ then we get a $C^M$-curve, since the continuous linear functionals on $\mathbb R^\mathbb N$ are linear combinations of coordinate projections $\on{ev}_t$ with $t\in\mathbb N$. 
However, this curve cannot be $C_b^M$ as the argument above for $t>\rho$ shows.

\section{\label{nmb:2} Working up to cartesian closedness: More on non-quasianalytic functions}

In \cite{KMRc} we developed convenient calculus for $C^M$ where
$M$  was  log-convex,  increasing,  derivation  closed,  and of
moderate growth for the exponential law.
In   this   paper   we   describe   quasianalytic  mappings  as
intersections of non-quasianalytic classes $C^L$, but we cannot
assume that $L$ is derivation closed.
Thus  we  need stronger versions of many results of \cite{KMRc}
for  non-quasianalytic  $L$  which  are  {\bf  not}  derivation
closed, and sometimes even {\bf not} log-convex.
This  section  collects  an almost minimal set of results which
allow  to  prove cartesian closedness for certain quasianalytic
function classes.

\begin{lemma}[cf.\   {\cite[3.3]{KMRc}}]   \label{nmb:2.1}  Let
$M=(M_k)_{k\in\mathbb N}$ be a sequence of positive numbers
and let $E$ be a convenient vector space such that there exists a Baire vector
space topology on the dual $E^*$ for which the point evaluations $\on{ev}_x$
are continuous for all $x\in E$.
Then a curve $c:\mathbb R\to E$ is $C^M$ if and
only if $c$ is $C_b^M$.
\end{lemma}

\begin{demo}{Proof}
Let $K$ be compact in $\mathbb R$ and $c$ be a $C^M$-curve.
We consider the sets
\[
A_{\rho,C} :=\Bigl\{\ell\in E^*: \frac{|\ell(c^{(k)}(x))|}{\rho^{k}\,
k!\, M_{k}}\le C\text{ for all }k\in \mathbb{N}, x\in K\Bigr\}
\]
which are closed subsets in $E^*$ for the given Baire topology. 
We have
$\bigcup_{\rho,C}A_{\rho,C}= E^*$. 
By the Baire property there exists $\rho$
and $C$ such that the interior $U$ of $A_{\rho,C}$ is non-empty. 
If 
$\ell_0\in U$ then for each $\ell\in E^*$ there is a $\delta>0$ such that 
$\delta\ell\in U-\ell_0$ and hence for all $x\in K$ and all $k$ we have
\begin{align*}
|(\ell\circ c)^{(k)}(x)| \le \tfrac1\delta \left(|((\delta\,\ell+\ell_0)\circ c)^{(k)}(x)| +
|(\ell_0\circ c)^{(k)}(x)|\right) \le \tfrac{2C}{\delta}\,\rho^{k}\,k!\,
M_{k}.
\end{align*}
So the set \[
\left\{\frac{c^{(k)}(x)}{\rho^{k}\, k!\, M_{k}}: k\in\mathbb{N}, x\in K
\right\}
\]
is weakly bounded in $E$ and hence bounded. 
\qed\end{demo}

\begin{lemma}[cf.\  {\cite[3.4]{KMRc}}]\label{nmb:2.2}%
Let $M=(M_k)_{k\in\mathbb N}$ be a sequence of positive numbers
and let $E$ be a Banach space.
For a smooth curve $c:\mathbb R\to E$ the following are equivalent.
\begin{enumerate}
\item[(1)] $c$ is $C^M=C^M_b$.
\item[(2)] For each sequence $(r_k)$ with $r_k\,\rho^k\to
0$ for all $\rho>0$, and each compact set $K$ in $\mathbb R$, the set 
$\{\frac1{k!M_k}\,c^{(k)}(a)\,r_k: a\in K, k\in \mathbb{N}\}$ is bounded in $E$.
\item[(3)] For each sequence $(r_k)$ satisfying $r_k>0$, $r_kr_\ell\geq r_{k+\ell}$, 
and $r_k\,\rho^k\to 0$ for all $\rho>0$, and each compact set $K$ in $\mathbb R$, 
there exists an
$\delta>0$ such that $\{\frac1{k!M_k}\,c^{(k)}(a)\,r_k\,\delta^k: a\in K, k\in \mathbb{N}\}$ 
is bounded in $E$. 
\end{enumerate}
\end{lemma}

\demo{Proof}
\thetag{1} $\implies$ \thetag{2}
For $K$, there exists $\rho>0$ such that
\[
\left\|\frac{c^{(k)}(a)}{k!\,M_k}r_k \right\|_E 
= \left\|\frac{c^{(k)}(a)}{k!\,\rho^k\,M_k}\right\|_E\cdot |r_k\rho^k| 
\]
is bounded uniformly in $k\in \mathbb N$ and $a\in K$ by \thetag{\ref{nmb:2.1}}. 
\thetag{2} $\implies$ \thetag{3} Use $\delta=1$.

\thetag{3} $\implies$ \thetag{1} Let $a_k:=\sup_{a\in K}\|\frac1{k!\,M_k}\,c^{(k)}(a)\|_E$. 
Using (4$\Rightarrow$1) in \cite[9.2]{KM97} these are the coefficients 
of a power series with positive radius of convergence.
Thus $a_k/\rho^k$ is bounded for some $\rho>0$. 
\qed\enddemo

\begin{lemma}[cf.\  {\cite[3.5]{KMRc}}]\label{nmb:2.3}
Let $M=(M_k)_{k\in\mathbb N}$ be a sequence of positive numbers.
Let $E$ be a convenient vector space, and let $\mathcal S$ be a
family of bounded linear functionals on $E$ which together detect bounded
sets (i.e., $B\subseteq E$ is bounded if and only if $\ell(B)$ is bounded for all
$\ell\in\mathcal S$). 
Then a curve $c:\mathbb R\to E$ is $C^M$ if and only if $\ell\circ c:\mathbb R\to \mathbb R$ is $C^M$
for all $\ell\in\mathcal S$.
\end{lemma}

\begin{demo}{Proof}
For smooth curves this follows from \cite[2.1, 2.11]{KM97}.
By \thetag{\ref{nmb:2.2}}, for $\ell\in \mathcal S$, 
the function $\ell\circ c$ is $C^M$ if and only if: 
\begin{enumerate}
\item
For each sequence $(r_k)$ with $r_k\,t^k\to
0$ for all $t>0$, and each compact set $K$ in $\mathbb R$, the set 
$\{\frac1{k!M_k}\,(\ell\circ c)^{(k)}(a)\,r_k: a\in K, k\in \mathbb{N}\}$ is bounded. 
\end{enumerate}
By (1) the curve $c$ is $C^M$ if and only if the set 
$\{\frac1{k!M_k}\,c^{(k)}(a)\,r_k: a\in K, k\in \mathbb{N}\}$ is
bounded in $E$. 
By (1) again this is in turn equivalent to $\ell\circ c\in C^M$ for all 
$\ell\in\mathcal S$, since $\mathcal S$ detects bounded sets.
\qed\end{demo}

\begin{corollary}\label{nmb:2.4}
Let  $M=(M_k)_{k\in\mathbb  N}$  be  a non-quasianalytic weight
sequence  or an $\mathcal L$-intersectable quasianalytic weight
sequence.
Let $U$ be $c^\infty$-open in a convenient vector space $E$, and let 
$\mathcal S=\{\ell:F\to F_{\ell}\}$ be a
family of bounded linear mappings between convenient vector spaces which together detect bounded
sets.
Then a mapping $f:U\to F$ is $C^M$ if and only if $\ell\circ f$ is $C^M$
for all $\ell\in\mathcal S$.

In particular, a mapping $f:U\to L(G,H)$ is $C^M$ if and only if 
$\on{ev}_v\circ f:U\to H$ is $C^M$ for each $v\in G$, where $G$ and $H$ are convenient vector spaces. 
\end{corollary}

This result is not valid for $C_b^M$ instead of $C^M$, by a variant of \thetag{\ref{nmb:1.11}}: Replace $C^M(\mathbb R,\mathbb R)$ by $\mathbb R^{\mathbb N}$.

\begin{demo}{Proof} First, let $M$ be non-quasianalytic.
By composing with curves we may reduce to $U=E=\mathbb R$.
By  composing  each $\ell\in\mathcal S$ with all bounded linear
functionals  on  $F_\ell$  we  get  a  family of bounded linear
functionals     on     $F$    to    which    we    can    apply
\thetag{\ref{nmb:2.3}}.
For quasianalytic $M$ the result follows by definition.
The  case  $F=L(G,H)$  follows  since  the $\on{ev}_v$ together
detect  bounded  sets,  by  the  uniform  boundedness principle
\cite[5.18]{KM97}.
\qed\end{demo}

\subsection{\label{nmb:2.5}$C^L$-curve lemma {\rm (cf.\ \cite[3.6]{KMRc})}}
A sequence $x_n$ in a locally convex space $E$ %
is said to be
{\it Mackey convergent} to $x$, if there exists some $\lambda_n\nearrow \infty$
such that $\lambda_n(x_n-x)$ is bounded. 
If we fix $\lambda=(\lambda_n)$ we say that $x_n$ is $\lambda$-converging. 

\begin{lemma*} Let $L$ be a non-quasianalytic weight sequence.
Then there exist sequences $\lambda_k\to 0$, $t_k\to t_\infty$, $s_k>0$ in 
$\mathbb R$ with the following property:

For $1/\lambda=(1/\lambda_n)$-converging sequences $x_n$ and $v_n$
in a convenient vector space $E$ there exists a strong uniform $C^L$-curve $c:\mathbb R\to E$ with
$c(t_k+t)= x_k+t.v_k$ for $|t|\leq s_k$.
\end{lemma*}

\begin{demo}{Proof}
Since $C^L$ is not quasianalytic we have $\sum_k 1/(k!L_k)^{1/k}<\infty$ by \thetag{\ref{nmb:1.2}}. 
We choose
another non-quasianalytic weight sequence $\bar L= (\bar L_k)$ with $(L_k/\bar L_k)^{1/k}\to \infty$. 
By   \cite[2.3]{KMRc}   there   is   a   $C^{\bar  L}$-function
$\phi:\mathbb R\to [0,1]$ which is 0 on $\{t: |t|\ge \frac12\}$
and which is 1 on
$\{t: |t|\le \frac13\}$, i.e.\  there exist $\bar C,\rho>0$ such that
\[
|\phi^{(k)}(t)|\leq \bar C\,\rho^k\,k!\,\bar L_k
\quad\text{ for all }t\in\mathbb R\text{ and }k\in\mathbb{N}.
\]
For $x,v$ in an absolutely convex bounded set $B\subseteq E$ and $0<T\le 1$ the curve
$c:t\mapsto \phi(t/T)\cdot (x+t\,v)$ satisfies
(cf.\ \cite[Lemma 2]{Boman67}):
\begin{align*}
c^{(k)}(t)&= T^{-k}\phi^{(k)}(\tfrac{t}{T}).(x+t.v)+ k\,T^{1-k}\,\phi^{(k-1)}(\tfrac{t}{T}).v
\\&
\in T^{-k}\bar C\,\rho^k\,k!\,\bar L_k(1+\tfrac{T}{2}).B+
k\,T^{1-k}\,\bar C\,\rho^{k-1}\,(k-1)!\,\bar L_{k-1}.B
\\&
\subseteq T^{-k}\bar C\,\rho^k\,k!\,\bar L_k(1+\tfrac{T}{2}).B+
T\,T^{-k}\,\bar C\tfrac1\rho\,\rho^{k}\,k!\,\bar L_{k}.B
\\&
\subseteq \bar C(\tfrac32+ \tfrac 1\rho)\,T^{-k}\,\rho^{k}\,k!\,\bar L_{k}.B
\end{align*}
So there are $\rho,C:=\bar C(\tfrac32+ \tfrac 1\rho)>0$ which do not depend on $x,v$ and $T$ such that
$c^{(k)}(t)\in C\, T^{-k}\,\rho^k\,k!\,\bar L_k. B$ for all $k$ and $t$.

Let $0<T_j\leq 1$ with $\sum_j T_j<\infty$
and $t_k:=2\sum_{j<k}T_j+T_k$.
We choose the $\lambda_j$ such that
$0<\lambda_j/T_j^k\leq     L_k/\bar     L_k$     (note     that
$T_j^k\,L_k/\bar  L_k\to\infty$  for  $k\to\infty$) for all $j$
and $k$, and that
$\lambda_j/T_j^k\to 0$ for $j\to\infty$ and each $k$.

Without loss we may assume that $x_n\to 0$. 
By  assumption  there exists a closed bounded absolutely convex
subset $B$ in $E$ such that $x_n,v_n\in \lambda_n\cdot B$.
We consider
$c_j:t\mapsto \phi\bigl((t-t_j)/T_j\bigr)\cdot \bigl(x_j+(t-t_j)\,v_j\bigr)$
and $c:=\sum_j c_j$.
The $c_j$ have disjoint support $\subseteq [t_j-T_j,t_j+T_j]$, hence $c$ is 
$C^\infty$ on $\mathbb R\setminus\{t_\infty\}$ with
\[
c^{(k)}(t) \in C\, T_j^{-k} \,\rho^{k} k!\bar L_k\,\lambda_j\cdot B \quad\text{ for }|t-t_j|\leq T_j.
\]
Then
\[
\|c^{(k)}(t)\|_{B}\le C\,\rho^{k}\,k!\bar L_k\,\frac{\lambda_j}{T_j^{k}} 
\leq C \rho^{k} k!\bar L_k\,\frac{L_k}{\bar L_k} = C\,\rho^{k}\,k! L_k
\]
for $t\ne t_\infty$.
Hence $c:\mathbb R\to E_B$ %
is smooth at $t_\infty$ as well, and is strongly $C^L$ by
the following lemma.%
\qed\end{demo}

\begin{lemma}[cf.\  {\cite[3.7]{KMRc}}]\label{nmb:2.6}
Let $c:\mathbb R\setminus\{0\}\to E$ be strongly $C^L$ in the sense that $c$ is smooth and for all
bounded $K \subset \mathbb R\setminus \{0\}$ there exists $\rho>0$ such that
\[
\left\{\frac{c^{(k)}(x)}{ \rho^{k} \, k! \, L_{k}}: k \in \mathbb{N}, x \in K\right\}\text{ is bounded in }E.
\]
Then $c$ has a unique extension to a strongly $C^L$-curve on $\mathbb R$.
\end{lemma}

\begin{demo}{Proof}
The curve $c$ has a unique extension to a smooth curve by \cite[2.9]{KM97}.
The strong $C^L$ condition extends by continuity. 
\qed\end{demo}

\begin{theorem}[cf.\  {\cite[3.9]{KMRc}}]\label{nmb:2.7}
Let $L=(L_k)$ be a non-quasianalytic weight sequence.
Let $U\subseteq E$ be $c^\infty$-open in a convenient vector space, let $F$
be a Banach space and $f:U\to F$ a mapping. 
Furthermore,
let $\overline{L}\le L$ be another non-quasianalytic weight sequence. 
Then the following statements are equivalent:
\begin{enumerate}
\item[(1)] $f$ is $C^L$, i.e.\ $f\circ c$ is $C^L$ for all $C^L$-curves $c$.
\item[(2)] $f|_{U\cap E_B}:E_B\supseteq U\cap E_B\to F$ is $C^L$ for each closed bounded absolutely convex $B$ in $E$.
\item[(3)] $f\circ c$ is $C^L$ for all $C_b^{\overline{L}}$-curves $c$.
\item[(4)] $f\in C_b^L(U,F)$.
\end{enumerate}
\end{theorem}

\begin{demo}{Proof}
\thetag{1} $\implies$ \thetag{2} is clear, since $E_B\to E$ is continuous and linear, 
hence all $C^L$-curves $c$ into the Banach space $E_B$ are also $C^L$ into $E$ 
and hence $f\circ c$ is $C^L$ by assumption.

\thetag{2} $\implies$ \thetag{3} is clear, since $C_b^{\overline{L}}\subseteq C^L$.

\thetag{3} $\implies$ \thetag{4}
Without loss let $E=E_B$ be a Banach space. 
For each $v\in E$ and $x\in U$ the iterated directional derivative $d_v^kf(x)$ 
exists since $f$ is $C^L$ along affine lines. 
To show that $f$ is smooth it suffices to check that
$d_{v_n}^kf(x_n)$ is bounded for each $k\in \mathbb{N}$ and each Mackey
convergent sequences $x_n$ and $v_n\to 0$, by \cite[5.20]{KM97}. 
For contradiction let us assume that there exist $k$ and sequences $x_n$
and $v_n$ with $\|d_{v_n}^kf(x_n)\|\to \infty$. 
By passing to a subsequence we may assume that $x_n$ and $v_n$ are $(1/\lambda_n)$-converging
for the $\lambda_n$ from \thetag{\ref{nmb:2.5}} for the weight sequence $\overline{L}$.
Hence there exists a $C_b^{\overline{L}}$-curve $c$
in $E$ and with $c(t+t_n)=x_{n}+t.v_{n}$ for $t$ near 0 for each
$n$ separately, and for $t_n$ from \thetag{\ref{nmb:2.5}}. 
But then $\|(f\circ c)^{(k)}(t_n)\|=\|d_{v_{n}}^k f(x_{n})\|\to \infty$, a contradiction.
So $f$ is smooth. 
Assume for contradiction that the boundedness condition in
\thetag{4} does not hold:
There exists a compact set $K\subseteq U$ such that for each $n\in\mathbb{N}$ there are 
$k_n\in\mathbb{N}$, $x_n\in K$, and $v_n$ with
$\|v_n\|=1$ such that
\[
\|d_{v_n}^{k_n}f(x_n)\|>{k_n}!\,L_{k_n}\,\left(\frac1{\lambda_n^2}\right)^{k_n+1},
\]
where we used $C=\rho:=1/\lambda_n^2$ with the $\lambda_n$ from \thetag{\ref{nmb:2.5}} 
for the weight sequence $\overline{L}$.
By passing to a subsequence (again denoted $n$) we may assume that the $x_n$ are
$1/\lambda$-converging, thus
there exists a $C_b^{\overline{L}}$-curve $c:\mathbb R\to E$ with $c(t_n+t)=x_n+t.\lambda_n.v_n$ 
for $t$ near 0 by \thetag{\ref{nmb:2.5}}.
Since
\[
(f\circ c)^{(k)}(t_n)= \lambda_n^k d_{v_n}^k f(x_n),
\]
we get \[
\left(\frac{\|(f\circ c)^{(k_n)}(t_n)\|}{{k_n}!L_{k_n}}\right)^{\frac1{k_n+1}}
=\left(\lambda_n^{k_n} \frac{\|d_{v_n}^{k_n}f(x_n)\|}{{k_n}!L_{k_n}}\right)^{\frac1{k_n+1}}
> \frac1{\lambda_n^{\frac{k_n+2}{k_n+1}}} \to\infty,
\]
a contradiction to $f\circ c\in C^L$.

\thetag{4} $\implies$ \thetag{1}
We have to show that $f\circ c$ is $C^L$ for each $C^L$-curve $c:\mathbb R\to
E$. 
By \thetag{\ref{nmb:2.2}.3} it suffices to show that for each sequence $(r_k)$ 
satisfying $r_k>0$, $r_kr_\ell\geq r_{k+\ell}$, and $r_k\,t^k\to 0$ for all $t>0$, 
and each compact interval $I$ in $\mathbb R$, there exists an
$\epsilon>0$ such that 
$\{\frac1{k!L_k}\,(f\circ c)^{(k)}(a)\,r_k\,\epsilon^k: a\in I, k\in \mathbb{N}\}$ is bounded.

By \thetag{\ref{nmb:2.2}.2} applied to $r_k2^k$ instead of $r_k$, for each 
$\ell\in E^*$, each sequence $(r_k)$ with $r_k\,t^k\to 0$ for all $t>0$, 
and each compact interval $I$ in $\mathbb R$
the set $\{\frac1{k!L_k}\,(\ell\circ c)^{(k)}(a)\,r_k\,2^k: a\in I, k\in \mathbb{N}\}$ 
is bounded in $\mathbb R$.
Thus $\{\frac1{k!L_k}\,c^{(k)}(a)\,r_k\,2^k: a\in I, k\in \mathbb{N}\}$ 
is contained in some closed absolutely convex $B\subseteq E$. 
Consequently, $c^{(k)}:I\to E_B$ is smooth and hence
$K_k:=\{\frac1{k!L_k}\,c^{(k)}(a)\,r_k\,2^k: a\in I\}$ is
compact in $E_B$ for each $k$. 
Then each sequence $(x_n)$ in the set
\[
K:=\left\{\frac1{k!L_k}\,c^{(k)}(a)\,r_k: a\in I, k\in \mathbb{N}\right\}
=\bigcup_{k\in\mathbb{N}}\frac1{2^k} K_k
\]
has a cluster point in $K\cup\{0\}$: either there is a subsequence in one
$K_k$, or $2^{k_n}x_{k_n}\in K_{k_n}\subseteq B$ for $k_n\to\infty$, hence
$x_{k_n}\to 0$ in $E_B$. 
So $K\cup\{0\}$ is compact. 
By Fa\`a di Bruno (\cite{FaadiBruno1855} for the 1-dimensional version, $k\ge 1$)
\begin{align*}
\frac{(f\circ c)^{(k)}(a)}{k!} = \sum_{j\ge 1} \sum_{\substack{\alpha\in \mathbb{N}_{>0}^j\\ \alpha_1+\dots+\alpha_j =k}}
\frac{1}{j!}d^jf(c(a))\Big( \frac{c^{(\alpha_1)}(a)}{\alpha_1!},\dots,
\frac{c^{(\alpha_j)}(a)}{\alpha_j!}\Big)
\end{align*}
and \thetag{\ref{nmb:1.1}.2} for $a\in I$ and $k\in \mathbb{N}_{>0}$ we have
\begin{align*}
&\left\|\frac1{k!L_k}\,(f\circ c)^{(k)}(a)\,r_k \right\|\le
\\&
\le \sum_{j\ge 1} L_1^j \sum_{\substack{\alpha\in \mathbb{N}_{>0}^j\\ \alpha_1+\dots+\alpha_j =k}}
\frac{\|d^jf(c(a))\|_{L^j(E_B,F)}}{j!L_j}\prod_{i=1}^j
\frac{\|c^{(\alpha_i)}(a)\|_B \, r_{\alpha_i}}{\alpha_i!L_{\alpha_i}}
\\&
\le \sum_{j\ge 1} L_1^j \binom{k-1}{j-1} C\,\rho^j\, \frac1{2^k}
= L_1 \rho(1+L_1\,\rho)^{k-1} C\, \frac1{2^k}.
\end{align*}
So $
\left\{\frac1{k! L_k}\,(f\circ c)^{(k)}(a)\,\Big(\frac{2}{1+L_1\,\rho}\Big)^k\,r_k: a\in I, k\in \mathbb{N}\right\}
$
is bounded as required.
\qed\end{demo}

\begin{corollary}\label{nmb:2.8}
Let $L=(L_k)$ be a non-quasianalytic weight sequence.
Let $U\subseteq E$ be $c^\infty$-open in a convenient vector space, let $F$
be a convenient vector space and $f:U\to F$ a mapping. 
Furthermore,
let $\overline{L}\le L$ be a non-quasianalytic weight sequence.
Then the following statements are equivalent:
\begin{enumerate}
\item[(1)] $f$ is $C^L$.
\item[(2)] $f|_{U\cap E_B}:E_B\supseteq U\cap E_B\to F$ is $C^L$ 
for each closed bounded absolutely convex $B$ in $E$.
\item[(3)] $f\circ c$ is $C^L$ for all $C_b^{\overline{L}}$-curves $c$.
\item[(4)] $\pi_V\circ f\in C_b^L(U,\mathbb R)$ for each absolutely convex 0-neighborhood 
$V\subseteq F$, where $\pi_V:F\to F_V$ denotes the natural mapping.
\end{enumerate}
\end{corollary}

\begin{demo}{Proof}
Each of the statements holds for $f$ if and only if it holds for $\pi_V\circ f$ for
each absolutely convex 0-neighborhood $V\subseteq F$. 
So the corollary follows from \thetag{\ref{nmb:2.7}}.
\qed\end{demo}

\begin{theorem}[Uniform boundedness principle for $C^M$, cf.\ {\cite[4.1]{KMRc}}]\label{nmb:2.9}
Let $M=(M_k)$ be a non-quasianalytic weight sequence or an 
$\mathcal L$-intersectable quasianalytic weight sequence. 
Let $E$, $F$, $G$ be convenient vector spaces and let $U\subseteq F$ be
$c^\infty$-open. 
A linear mapping $T:E\to C^M(U,G)$ is bounded if and only
if $\on{ev}_x\circ T: E\to G$ is bounded for every $x\in U$. \end{theorem}

\begin{demo}{Proof} Let first $M$ be non-quasianalytic.
For $x\in U$ and $\ell\in G^*$ the linear mapping $\ell\circ\on{ev}_x = C^M(x,\ell):C^M(U,G)\to
\mathbb R$ is continuous, thus $\on{ev}_x$ is bounded. Therefore, 
if $T$ is bounded then so is $\on{ev}_x\circ T$.

Conversely, suppose that $\on{ev}_x\circ T$ is bounded for all $x\in U$. 
For each closed absolutely convex bounded $B\subseteq E$ we consider the Banach
space $E_B$. For each $\ell\in
G^*$, each $C^M$-curve $c:\mathbb R\to U$, each $t\in\mathbb R$, and each compact 
$K\subset \mathbb R$ the composite given by the following diagram
is bounded. \[
\xymatrix{
E\ar@{->}[0,1]^{T\qquad} &C^M(U,G)\ar@{->}[1,0]^{C^M(c,\ell)} \ar@{->}[0,2]
^{\on{ev}_{c(t)}} & &G\ar@{->}[1,0]^{\ell} \\
E_B\ar@{->}[-1,0] \ar@{->}[0,1] &C^M(\mathbb R,\mathbb R)\ar@{->}[0,1] &
\varinjlim _{\rho }C^M_\rho (K,\mathbb R)\ar@{->}[0,1]^{\qquad\on{ev}_t} &
\mathbb R \\
}
\]
By \cite[5.24, 5.25]{KM97} the map $T$ is bounded. In more detail:
Since $\varinjlim_{\rho}C^M_\rho(K,\mathbb R)$ is webbed, the
closed graph theorem \cite[52.10]{KM97} yields that the mapping 
$E_B\to \varinjlim_{\rho}C^M_\rho(K,\mathbb R)$ is continuous. Thus $T$ is
bounded.

For quasianalytic $M$ the result follows since the structure of a convenient vector space on 
$C^M(U,G)$ is the initial one with respect to all inclusions $C^M(U,G)\to C^L(U,G)$ for all 
$L\in \mathcal L(M)$. \qed\end{demo}

As a consequence we can show that the equivalences of \thetag{\ref{nmb:2.7}} and 
\thetag{\ref{nmb:2.8}} are not only valid for single
functions $f$ but also for the bornology of $C^M(U,F)$:

\begin{corollary}[cf.\ {\cite[4.6]{KMRc}}]\label{nmb:2.10}
Let $L=(L_k)$ be a non-quasianalytic weight sequence.
Let $E$ and $F$ be Banach spaces and let $U\subseteq E$ be open.
Then \begin{align*}
C^L(U,F) = C_b^L (U,F):= \varprojlim_K \varinjlim_\rho C^L_{K,\rho}(U,F) \end{align*}
as vector spaces with bornology.
Here $K$ runs through all compact subsets of $U$ ordered by inclusion and
$\rho$ runs through the positive real numbers.
\end{corollary}

\begin{demo}{Proof}
The second equality is by definition \thetag{\ref{nmb:1.10}}. 
The first equality, as vector spaces, is by \thetag{\ref{nmb:2.7}}.
By \thetag{\ref{nmb:1.10}} the space $C^L(U,F)$ is convenient.

The identity from right to left is continuous since $C^L(U,F)$ 
carries the initial structure with respect to the mappings
\[
C^L(c|_I,\ell):C^L(U,F)\to C^L(\mathbb R,\mathbb R)
=\varprojlim_{I\subseteq\mathbb R}\varinjlim_{\rho>0}C^L_{I,\rho}(\mathbb R,\mathbb R)
\to \varinjlim_{\rho>0}C^L_{I,\rho}(\mathbb R,\mathbb R),
\]
where $c$ runs through the 
$C^L=\mkern-16mu\overset{\text{\thetag{\ref{nmb:2.1}}}}{\cleaders\hbox{$\mkern-3mu\mathord=\mkern-3mu$}\hfill}\mkern-16mu=C_b^L$-curves, $\ell\in F^*$ and $I$ runs through the compact intervals in $\mathbb R$,
and for $K:=c(I)$ and $\rho':=(1+\rho\,\|c\|_{I,\sigma})\cdot\sigma$, where $\sigma>0$ is chosen such that
$\|c\|_{I,\sigma}<\infty$,
the mapping 
$C^L(c|_I,\ell):C^L_{K,\rho}(U,F)\to C^L_{I,\rho'}(\mathbb R,\mathbb R)
\to \varinjlim_{\rho'>0} C^L_{I,\rho'}(\mathbb R,\mathbb R)$ is continuous by \thetag{\ref{nmb:1.4}}.
These arguments are collected in the diagram:
\[\xymatrix@!C=1.6cm{
\varprojlim_IC^L_I(\mathbb R,\mathbb R)\ar@{={}}[0,1] \ar@{->}[1,0] &C^L(\mathbb R,\mathbb R) &
C^L(U,F)\ar@{->}[0,-1]^{C^L(c,\ell)} &C_b^L(U,F)\ar@{.>}[0,-1] \ar@{={}}[0,
1] &\varprojlim_K C^L_K(U,F)\ar@{->}[1,0] \\
C^L_I(\mathbb R,\mathbb R)\ar@{={}}[0,1] &\varinjlim_{\rho}C^L_{I,\rho}(\mathbb R,\mathbb R) & &
\varinjlim_{\rho}C^L_{K,\rho}(U,F)\ar@{.>}[0,-2] &C^L_K(U,F)\ar@{={}}[0,-1] \\
&C^L_{I,\rho'}(\mathbb R,\mathbb R)\ar@{->}[-1,0] & &C^L_{K,\rho}(U,F)\ar@{->}[-1,0] \ar@{->}[0,-2]_{C^L(c|_I,\ell)} & \\
}\]
The identity from left to right is bounded since the countable 
(take $\rho\in\mathbb N$) inductive limit $\varinjlim_\rho$ 
of the (non-Hausdorff) Banach spaces $C^L_{K,\rho}(U,F)$ is webbed and 
hence satisfies the $\mathcal{S}$-boundedness principle
\cite[5.24]{KM97} where $\mathcal{S}=\{\on{ev}_x:x\in U\}$, and by 
\cite[5.25]{KM97} the same is true for $C_b^L(U,F)$.
\qed\end{demo}

\begin{corollary}[cf.\ {\cite[4.4]{KMRc}}]\label{nmb:2.11}
Let $L=(L_k)$ be a non-quasianalytic weight sequence.
Let $E$ and $F$ be convenient vector spaces and let $U\subseteq E$ be $c^\infty$-open.
Then
\[
C^L(U,F)
=\varprojlim_{c\in C^L} C^L(\mathbb R,F)
=\varprojlim_{B\subseteq E} C^L(U\cap E_B,F)
=\varprojlim_{s\in C_b^L} C^L(\mathbb R,F)
\]
as vector spaces with bornology, where
$c$ runs through all $C^L$-curves in $U$,
$B$ runs through all bounded closed absolutely convex subsets of $E$,
and $s$ runs through all $C_b^L$-curves in $U$.
\end{corollary}
\begin{demo}{Proof}
The first and third inverse limit is formed with $g^*:C^L(\mathbb R,F)\to C^L(\mathbb R,F)$ for 
$g\in C^L(\mathbb R,\mathbb R)$ as connecting mappings.
Each element $(f_c)_c$ determines a unique function $f:U\to F$ given by 
$f(x):=(f\circ\on{const}_x)(0)$ with $f\circ c=f_c$ for all such curves $c$, 
and $f\in C^L$ if and only if $f_c\in C^L$ for all such $c$, by \thetag{\ref{nmb:2.8}}.
The second inverse limit is formed with $\on{incl}^*:C^L(U\cap E_B,F)\to C^L(U\cap E_{B'},F)$ for 
$B'\subseteq B$ as connecting mappings.
Each element $(f_B)_B$ determines a unique function $f:U\to F$ given by 
$f(x):=f_{[-1,1]x}(x)$ with $f|_{E_B}=f_B$ for all $B$, and $f\in C^L$ 
if and only if $f_B\in C^L$ for all such $B$, by \thetag{\ref{nmb:2.8}}.
Thus all equalities hold as vector spaces.

The first identity is continuous from left to right, since the family of 
$\ell_*:C^L(\mathbb R,F)\to C^L(\mathbb R,\mathbb R)$ with $\ell\in F^*$ is initial and 
$C^L(c,\ell)=\ell_*\circ c^*:C^L(U,F)\to C^L(\mathbb R,\mathbb R)$ is continuous and linear by definition.

Continuity for the second one from left to right is obvious, since $C^L$-curves in $U\cap E_B$ 
are $C^L$ into $U\subseteq E$.

In order to show the continuity of the last identity from left to right choose a $C_b^L$-curve $s$
in $U$, an $\ell\in F^*$ and a compact interval $I\subseteq\mathbb R$. 
Then there exists a bounded absolutely convex closed $B\subseteq E$ such that $s|_I$ is $C_b^L=C^L$ into $U\cap E_B$, 
hence
$C^L(s|_I,\ell):C^L(U,F)\to C^L(I,\mathbb R)$ factors by \thetag{\ref{nmb:1.4}} as continuous linear mapping 
$(s|_I)^*:C_b^L(U\cap E_B,\mathbb R)\to C^L(I,\mathbb R)$ over
$C^L(U,F)\to C^L(U\cap E_B,F)\to C^L(U\cap E_B,\mathbb R)
=\mkern-16mu\overset{\text{\thetag{\ref{nmb:2.10}}}}{\cleaders\hbox{$\mkern-3mu\mathord=\mkern-3mu$}\hfill}\mkern-16mu
=C_b^L(U\cap E_B,\mathbb R)$. 
Since the structure of $C^L(\mathbb R,F)$ is initial with respect to 
$\on{incl}^*\circ\,\ell_*:C^L(\mathbb R,F)\to C^L(I,\mathbb R)$ the identity
$\varprojlim_{B\subseteq E} C^L(U\cap E_B,F)\to \varprojlim_{s\in C_b^L} C^L(\mathbb R,F)$ is continuous.

Conversely, the identity $\varprojlim_{s\in C_b^L} C^L(\mathbb R,F)\to C^L(U,F)$ is bounded, since
$C^L(\mathbb R,F)$ is convenient and hence also the inverse limit $\varprojlim_{s\in C_b^L} C^L(\mathbb R,F)$
and $C^L(U,F)$ satisfies the uniform boundedness theorem \thetag{\ref{nmb:2.9}} with respect to the point-evaluations $\on{ev}_x$ and they factor over $(\on{const}_x)^*:C^L(U,F)\to C^L(\mathbb R,F)$.
\qed\end{demo}

\section{\label{nmb:3} The exponential law for certain quasianalytic function classes}

We start with some preparations.
Let $Q=(Q_k)$ be an $\mathcal L$-intersectable quasianalytic weight sequence.
Let $E$ and $F$ be convenient vector spaces and let $U\subseteq E$ be $c^\infty$-open.

\begin{lemma}\label{nmb:3.1}
For Banach spaces $E$ and $F$
we have \[
C^Q(U,F)=C_b^Q(U,F)=\bigcap_{N\in \mathcal L_w(Q)}C_b^N(U,F)
\]
as vector spaces.
\end{lemma}

\begin{demo}{Proof} Since $Q$ is $\mathcal L$-intersectable
we have $\mathcal F^Q=\bigcap_{L\in\mathcal L(Q)}\mathcal F^L$.
Hence
\begin{align*}
C_b^Q(U,F) &=\{f\in C^\infty(U,F):
\forall K:(\sup_{x\in K}\|f^{(k)}(x)\|_{L^k(E,F)})_k\in\mathcal F^Q=\bigcap_{L\in\mathcal L(Q)}\mathcal F^L\} \\
&=\{f\in C^\infty(U,F):
\forall K\;\forall L\in\mathcal L(Q):(\sup_{x\in K}\|f^{(k)}(x)\|)_k\in\mathcal F^L\} \\
&=\{f\in C^\infty(U,F):
\forall L\in\mathcal L(Q)\;\forall K:(\sup_{x\in K}\|f^{(k)}(x)\|)_k\in\mathcal F^L\} \\
&= \bigcap_{L\in\mathcal L(Q)}C_b^L(U,F) =\mkern-16mu\overset{\text{\thetag{\ref{nmb:2.7}}}}{\cleaders\hbox{$\mkern-3mu\mathord=\mkern-3mu$}\hfill}\mkern-16mu= \bigcap_{L\in\mathcal L(Q)}C^L(U,F) = C^Q(U,F).
\\
C_b^Q(U,F) &=\mkern-16mu\overset{\text{\thetag{\ref{nmb:1.6}.1}}}{\cleaders\hbox{$\mkern-3mu\mathord=\mkern-3mu$}\hfill}\mkern-16mu=\bigcap_{L\in\mathcal L(Q)}C_b^L(U,F) \supseteq \bigcap_{L\in\mathcal L_w(Q)}C_b^L(U,F)
\supseteq C_b^Q(U,F). \qed
\end{align*}
\end{demo}

\begin{lemma}\label{nmb:3.2}
For log-convex non-quasianalytic
$L^1,L^2$ and weakly log-convex non-quasianalytic $N$ with
$N_{k+n} \le C^{k+n} L^1_k L^2_k$ for some positive constant $C$ and all $k,n \in \mathbb N$,
for Banach-spaces $E_1$ and $E_2$, and for $f\in C_b^N(U_1\times U_2,\mathbb R)$
we have
$f^\vee \in C^{L^1}(U_1,C_b^{L^2}(U_2,\mathbb R ))$.
\end{lemma}

\begin{demo}{Proof}
Since $f$ is $C_b^N$, by definition, for all compact $K_i\subseteq U_i$ there exists a $\rho>0$
such that for all $k,j\in\mathbb N$, $x_i\in K_i$ and $\|v_1\|=\dots=\|v_j\|=1=\|w_1\|=\dots=\|w_k\|$
we have
\begin{align*}
|\partial_2^k\partial_1^j f(x_1,x_2)&(v_1,\dots,v_j,w_1,\dots,w_k)|
\leq \rho^{k+j+1}(k+j)!N_{k+j} \\
&\leq \rho^{k+j+1}2^{k+j}k!j!C^{k+j}L^1_jL^2_k
= \rho (2C\rho)^j j! L^1_j \cdot (2C\rho)^k k! L^2_k.
\end{align*}
In particular $(\partial_1^jf)^\vee(K_1)(oE_1^k)$ is contained and bounded in 
$C_b^{L^2}(U_2,\mathbb R)$, where $oE_1$ denotes the unit ball in $E_1$, since
$d^k((\partial_1^jf)^\vee(x_1))(x_2)=\partial_2^k\partial_1^j f(x_1,x_2)$.

\noindent
{\bf Claim.} {\it If $f\in C_b^N$ then
$f^\vee:U_1\to C_b^{L^2}(U_2,\mathbb R)$ is $C^\infty$
with $d^jf^\vee=(\partial_1^jf)^\vee$.}
\newline
Since $C_b^{L^2}(U_2,\mathbb R)$ is a convenient vector space, by \cite[5.20]{KM97}
it is enough to show that the iterated unidirectional
derivatives $d^j_vf^\vee(x)$ exist, equal $\partial_1^jf(x,\quad)(v^j)$, 
and are separately bounded for $x$, resp.\  $v$, in compact subsets.
For $j=1$ and fixed $x$, $v$, and $y$ consider the smooth curve
$c:t\mapsto f(x+tv,y)$. By the fundamental theorem
\begin{align*}
\frac{f^\vee(x+tv)-f^\vee(x)}{t}(y)&-(\partial_1f)^\vee(x)(y)(v) = \frac{c(t)-c(0)}t-c'(0) \\
&= t\int_0^1 s \int_0^1 c''(tsr)\,dr\,ds \\
&= t\int_0^1 s\int_0^1 \partial_1^2f(x+tsrv,y)(v,v)\,dr \,ds.
\end{align*}
Since $(\partial_1^2f)^\vee(K_1)(oE_1^2)$ is bounded in $C_b^{L^2}(U_2,\mathbb R)$ for each compact subset
$K_1\subseteq U_1$ this expression is Mackey convergent %
to 0 in $C_b^{L^2}(U_2,\mathbb R)$, for $t\to 0$.
Thus $d_vf^\vee(x)$ exists and equals $\partial_1f(x,\quad)(v)$.

Now we proceed by induction, applying the same arguments as before to
$(d^j_vf^\vee)^\wedge : (x,y)\mapsto \partial_1^jf(x,y)(v^j)$ instead of $f$.
Again
$(\partial_1^2(d^j_vf^\vee)^\wedge)^\vee(K_1)(oE_1^2) = (\partial_1^{j+2}f)^\vee(K_1)(oE_1,oE_1,v,\dots,v)$
is bounded, and also the separated boundedness of $d^j_vf^\vee(x)$ follows.
So the claim is proved.

It remains to show that
$f^\vee:U_1\to C_b^{L^2}(U_2,\mathbb R)
:=\varprojlim_K\varinjlim_\rho C^{L^2}_{K,\rho}(U_2,\mathbb R)$ is $C^{L^1}$. 
By \thetag{\ref{nmb:2.4}}, it suffices to show that 
$f^\vee:U_1\to \varinjlim_\rho C^{L^2}_{K_2,\rho}(U_2,\mathbb R)$ is 
$C_b^{L^1}\subseteq C^{L^1}$ for all $K_2$, i.e., for all compact 
$K_2\subset U_2$ and $K_1\subset U_1$ there exists $\rho_1>0$ such that 
\[
\left\{\frac{d^kf^\vee(K_1)(v_1,\dots,v_k)}{k!\rho_1^{k}L^1_k}:k\in\mathbb N,\|v_i\|\leq 1\right\}
\text{ is bounded in }\varinjlim_\rho C^{L^2}_{K_2,\rho}(U_2,\mathbb R),
\]
or equivalently:
For all compact $K_2\subset U_2$ and $K_1\subset U_1$ there exist $\rho_1>0$ and $\rho_2>0$ such that 
\[
\left\{\frac{\partial_2^l\partial_1^k f(K_1,K_2)(v_1,\dots,v_{k+l})}{l!k!\rho_2^lL^2_l\rho_1^{k}L^1_k}:
k\in\mathbb N,l\in\mathbb N, \|v_i\|\leq 1\right\}
\text{ is bounded in }\mathbb R.
\]
For
$k_1\in\mathbb{N}$, $x\in K_1$, $\rho_i:=2C\rho$, and $\|v_i\|\le 1$
we get:
\begin{align*}
&\left\| \frac{d^{k_1}f^\vee (x)(v_1,\dots,v_{k_1})}{\rho_1^{k_1}\,k_1!\,L^1_{k_1}}\right\|_{K_2,\rho_2}= \\
&:=
\sup\Bigl\{\frac{|\partial_2^{k_2}\partial_1^{k_1}f(x,y)(v_1,\dots;w_1,\dots)|}{\rho_1^{k_1}\,k_1!\,L^1_{k_1}\,
\rho_2^{k_2}\,k_2!\,L^2_{k_2}}: k_2\in\mathbb{N}, y\in K_2,\|w_i\|\le 1
\Bigr\}
\\&
\le \sup\Bigl\{\frac{\frac{(k_1+k_2)!}{k_1!\,k_2!}C^{k_1+k_2} |\partial_2^{k_2}\partial_1^{k_1}f(x,y)(v_1,\dots;w_1,\dots)|}
{\rho_1^{k_1}\,\rho_2^{k_2}
\,(k_1+k_2)!
\,N_{k_1+k_2}}: k_2\in\mathbb{N}, y\in K_2,\|w_i\|\le 1\Bigr\}
\\&
\le \sup\Bigl\{\frac{(2C)^{k_1+k_2} |\partial^{(k_1,k_2)}f(x,y)(v_1,\dots;w_1,\dots)|}
{\rho_1^{k_1}\,\rho_2^{k_2}\,
\,(k_1+k_2)!\,
N_{k_1+k_2}}: k_2\in\mathbb{N}, y\in K_2,\|w_i\|\le 1\Bigr\}
\\&
=
\sup\Bigl\{\frac{|\partial^{(k_1,k_2)}f(x,y)(v_1,\dots;w_1,\dots)|}
{\rho^{k_1+k_2}\,
(k_1+k_2)!\,N_{k_1+k_2}}: k_2\in\mathbb{N}, y\in K_2:\|w_i\|\leq 1\Bigr\}
\le \rho
\end{align*}
So $f^\vee$ is $C^{L_1}$.
\qed\end{demo}

\begin{theorem}[Cartesian closedness]\label{nmb:3.3} 
Let $Q=(Q_k)$ be an $\mathcal L$-intersectable quasianalytic weight sequence of moderate growth.
Then the
category of $C^Q$-mappings between convenient real vector spaces is cartesian closed.
More precisely, for convenient vector spaces $E_1$, $E_2$ and $F$
and $c^\infty$-open sets $U_1\subseteq E_1$ and $U_2\subseteq E_2$ a mapping 
$f:U_1\times U_2\to F$ is $C^Q$ if and only
if $f^\vee :U_1\to C^Q(U_2,F)$ is $C^Q$.
\end{theorem}

Actually, we prove that the direction ($\Leftarrow$) holds without the assumption of moderate growth.

\begin{demo}{Proof} ($\Rightarrow$)
Let $f:U_1\times U_2\to F$ be $C^Q$, i.e.\ $C^L$ for all $L\in\mathcal L(Q)$.
Since $(E_i)_{B_i}\to E_i$ is bounded and linear and since $C^L$ is closed under composition we get that
$\ell\circ f:(U_1\cap (E_1)_{B_1})\times(U_2\cap (E_2)_{B_2})\to \mathbb R$ is 
$C^L=C_b^L$ (by \thetag{\ref{nmb:2.7}} since $(E_i)_{B_i}$ are Banach-spaces)
for $\ell\in F^*$, arbitrary bounded closed
$B_i\subseteq E_i$ and all $L\in\mathcal L(Q)$.
Hence $\ell\circ f$ is $C_b^L$ even for all $L\in\mathcal L_w(Q)$ by \thetag{\ref{nmb:3.1}}.
For arbitrary $L^1,L^2\in \mathcal L(Q)$, by \thetag{\ref{nmb:1.6}.3} and \thetag{\ref{nmb:1.6}.2}, 
there exists an $N\in\mathcal L_w(Q)$ with
$N_{k+n} \le C^{k+n} L^1_k L^2_n$ for some positive constant $C$ and all $k,n \in \mathbb N$.
Thus $\ell\circ f:(U_1\cap (E_1)_{B_1})\times(U_2\cap (E_2)_{B_2})\to \mathbb R$ is $C_b^N$.
By \thetag{\ref{nmb:3.2}}, the function
$(\ell\circ f)^\vee :U_1\cap (E_1)_{B_1}\to C_b^{L^2}(U_2\cap (E_2)_{B_2},\mathbb R)$
is $C^{L^1}$.
Since the cone 
\[
C^Q(U_2,F)\to C^{L^2}(U_2,F)-\raisebox{0.1pt}{$\mkern-16mu\frac{\;\;C^{L^2}(i_{B_2}, \ell)\;}{\;\;\;}\mkern-16mu$}\to C^{L^2}(U_2\cap (E_2)_{B_2},\mathbb R)
=C_b^{L^2}(U_2\cap (E_2)_{B_2},\mathbb R),
\]
with
$L_2\in\mathcal L(Q)$, $\ell\in F^*$, and bounded closed $B_2\subseteq E_2$, generates the bornology
by \thetag{\ref{nmb:2.11}}, and since obviously $f^\vee (x)=f(x,\quad)\in C^Q(U_2,F)$,
we have that
$f^\vee :U_1\cap (E_1)_{B_1}\to C^Q(U_2,F)$
is $C^{L^1}$, by \thetag{\ref{nmb:2.4}}. From this we get by \thetag{\ref{nmb:2.8}} that
$f^\vee :U_1\to C^Q(U_2,F)$ is $C^{L^1}$ for all $L^1\in\mathcal L(Q)$, i.e., $f^\vee :U_1\to C^Q(U_2,F)$ is $C^Q$ as required. The whole argument above is collected in the following diagram where $U^i_{{B_i}}$ stands for $U_i\cap E_{B_i}$:
\[
\xymatrix{
U^1\times U^2\ar@{->}[0,2]^(0.6){f\in C^Q} &&F\ar@{->}[1,0]_{\ell} &U^1\ar@{->}[0,1]^(0.4){f^\vee\in C^{L^1}}_(0.3){\text{\thetag{\ref{nmb:2.8}}}} &C^Q(U^2,F) %
\ar@{->}[0,1] & C^{L^2}(U^2,F) \ar@{->}[1,0]_{\ell_*\circ\,\on{incl}_2^*}^{\text{\thetag{\ref{nmb:2.11}}}} \\
U^1_{{B_1}}\times U^2_{{B_2}}\ar@{->}[-1,0]_{\on{incl}} \ar@{->}[0,2]^(0.6){f\in C^Q\subseteq C_b^N}_(0.6){\text{\thetag{\ref{nmb:3.1}}}} && \mathbb R \ar@{}[0,1]|{\Longrightarrow} &U^1_{{B_1}}\ar@{->}[-1,0]_{\on{incl}_1} \ar@{->}[0,1]^(0.4){C^{L^1}}_(0.4){\text{\thetag{\ref{nmb:3.2}}}} \ar@{->}[-1,1]_(0.7){\text{\thetag{\ref{nmb:2.3}}}}
& C^{L^2}(U^2_{{B_2}},\mathbb R)\ar@{={}}[0,1] & C_b^{L^2}(U^2_{{B_2}},\mathbb R)
}
\]

\medskip\noindent
($\Leftarrow$)
Let, conversely, $f^\vee :U_1\to C^Q(U_2,F)$ be $C^Q$, i.e.,
$C^{L}$ for all $L\in \mathcal L(Q)$. By the description of the structure of $C^Q(U,F)$ 
in \thetag{\ref{nmb:1.10}} the mapping $f^\vee :U_1\to C^L(U_2,F)$ is $C^L$. 
We now conclude that $f:U_1\times U_2\to F$ is $C^L$; this direction of cartesian closedness for
$C^L$ holds even if $L$ is not of moderate growth, see \cite[5.3]{KMRc} and its proof.
This is true for all $L\in\mathcal L(Q)$. Hence $f$ is $C^Q$. \qed\end{demo}

\begin{corollary}\label{nmb:3.4}
Let $Q$ be an $\mathcal L$-intersectable quasianalytic weight sequence of
moderate growth.
Let $E$, $F$, etc., be convenient vector spaces and let $U$ and $V$ be
$c^\infty$-open subsets of such. Then we have:
\newline \thetag{1} The exponential law holds: \[
C^Q(U,C^Q(V,G)) \cong C^Q(U\times V, G)
\]
\indent\indent is a linear $C^Q$-diffeomorphism of convenient vector spaces. \newline
The following canonical mappings are $C^Q$.
\begin{align*}
&\operatorname{ev}: C^Q(U,F)\times U\to F,\quad \operatorname{ev}(f,x) = f(x)
\tag{2}\\&
\operatorname{ins}: E\to C^Q(F,E\times F),\quad
\operatorname{ins}(x)(y) = (x,y)
\tag{3}\\&
(\quad)^\wedge :C^Q(U,C^Q(V,G))\to C^Q(U\times V,G)
\tag{4}\\&
(\quad)^\vee :C^Q(U\times V,G)\to C^Q(U,C^Q(V,G))
\tag{5}\\&
\operatorname{comp}:C^Q(F,G)\times C^Q(U,F)\to C^Q(U,G)
\tag{6}\\&
C^Q(\quad,\quad):C^Q(F,F_1)\times C^Q(E_1,E)\to C^Q\Bigl(C^Q(E,F),C^Q(E_1,F_1)\Bigr)
\tag{7}\\&
\qquad (f,g)\mapsto(h\mapsto f\circ h\circ g)
\\&
\prod:\prod C^Q(E_i,F_i)\to C^Q\Bigl(\prod E_i,\prod F_i\Bigr)
\tag{8}\end{align*}
\end{corollary}

\begin{demo}{Proof}
This is a direct consequence of cartesian closedness \thetag{\ref{nmb:3.3}}. See \cite[5.5]{KMRc} or even \cite[3.13]{KM97} for the detailed arguments.
\qed\end{demo}

\section{\label{nmb:4} More on function spaces}

In  this  section we collect results for function classes $C^M$
where  $M$  is either a non-quasianalytic weight sequence or an
$\mathcal  L$-intersectable  quasianalytic  weight sequence. In
order to treat both cases simultaneously, the proofs will often
use  non-quasianalytic  weight  sequences  $L\ge  M$. These are
either  $M$  itself  if  $M$  is  non-quasianalytic  or  are in
$\mathcal   L(M)$   if   $M$   is   $\mathcal  L$-intersectable
quasianalytic.
In   both  cases  we  may  assume  without  loss  that  $L$  is
increasing, by \thetag{\ref{nmb:1.5}}.
\begin{proposition}\label{nmb:4.1}
Let  $M=(M_k)$  be  a  non-quasianalytic  weight sequence or an
$\mathcal  L$-intersectable quasianalytic weight sequence. Then
we have:
\begin{enumerate}
\item[(1)] Multilinear mappings between convenient vector spaces are $C^M$ if and only
if they are bounded.
\item[(2)]
If $f:E\supseteq U\to F$ is $C^M$, then the derivative $df:U\to
L(E,F)$ is $C^{M_{+1}}$, and also $(df)^\wedge :U\times E\to F$
is $C^{M_{+1}}$, where the space $L(E,F)$ of all bounded linear
mappings is considered with the topology of uniform convergence
on bounded sets.
\item[(3)]
The chain rule holds.
\end{enumerate}
\end{proposition}

\begin{demo}{Proof}
\thetag{1}
If $f$ is $C^M$ then it is smooth by \thetag{\ref{nmb:2.8}} and hence
bounded by \cite[5.5]{KM97}. Conversely, if $f$ is multilinear and bounded then
it  is  smooth, again by \cite[5.5]{KM97}. Furthermore, $f\circ
i_B$  is multilinear and continuous and all derivatives of high
order vanish. Thus condition \thetag{\ref{nmb:2.8}.4} is satisfied, so $f$ is $C^M$. 
\thetag{2} Since $f$ is smooth, by \cite[3.18]{KM97} the map 
$df:U\to L(E,F)$ exists and is smooth. Let $L\geq M_{+1}$ be a non-quasianalytic weight sequence and 
$c:\mathbb R\to U$ be a $C^L$-curve. 
We have to show that $t\mapsto df(c(t))\in L(E,F)$ is $C^L$. 
By the uniform boundedness principle \cite[5.18]{KM97} and by \thetag{\ref{nmb:2.3}} 
it suffices to show that the mapping $t\mapsto c(t)\mapsto \ell(df(c(t))(v))\in \mathbb R$ is $C^L$ for each
$\ell\in F^*$ and $v\in E$.
We are reduced to show that $x\mapsto \ell(df(x)(v))$ satisfies the
conditions of \thetag{\ref{nmb:2.7}}. By \thetag{\ref{nmb:2.7}} applied to $\ell\circ f$, 
for each $L\ge M$, each closed bounded absolutely convex $B$ in $E$,
and each $x\in U\cap E_B$ there are $r>0$, $\rho>0$, and $C>0$ such that
\[
\frac1{k!\,L_k}\|d^k(\ell\circ f\circ i_B)(a)\|_{L^k(E_B,\mathbb R)} \le C\,\rho^k
\]
for all $a\in U\cap E_B$ with $\|a-x\|_{B}\le r$ and all $k\in\mathbb{N}$. 
For $v\in E$ and those $B$ containing $v$ we then have:
\begin{align*}
\|d^k(d(\ell\circ f)(&\quad)(v))\circ i_B)(a)\|_{L^k(E_B,\mathbb R)} 
=\|d^{k+1}(\ell\circ f\circ i_B)(a)(v,\dots)\|_{L^k(E_B,\mathbb R)} \\&
\le\|d^{k+1}(\ell\circ f\circ i_B)(a)\|_{L^{k+1}(E_B,\mathbb R)}\|v\|_{B} \le C\,\rho^{k+1}\, (k+1)!\, L_{k+1}
\\&
= C\rho\,((k+1)^{1/k}\rho)^{k}\, k!\, L_{k+1}
\leq C\rho\,(2\rho)^k\,k!\,(L_{+1})_k
\end{align*}
By \thetag{\ref{nmb:4.2}} below also $(df)^\wedge $is $C^{L_{+1}}$.

\thetag{3} This is valid even for all smooth $f$ by \cite[3.18]{KM97}. \qed\end{demo}

\begin{proposition}\label{nmb:4.2} Let $M=(M_k)$ be a non-quasianalytic weight sequence or an $\mathcal L$-intersectable quasianalytic weight sequence. \begin{enumerate}
\item[(1)]
For convenient vector spaces $E$ and $F$,
on $L(E,F)$ the following bornologies coincide which are induced by:
\begin{itemize}
\item The topology of uniform convergence on bounded subsets
of $E$.
\item The topology of pointwise convergence.
\item The embedding $L(E,F)\subset C^\infty(E,F)$.
\item The embedding $L(E,F)\subset C^M(E,F)$.
\end{itemize}
\item[(2)] Let $E$, $F$, $G$ be convenient vector spaces and let $U\subset E$ be
$c^\infty$-open. A mapping $f:U\times F\to G$ which is linear in
the second variable is $C^M$ if and only if $f^\vee:U\to L(F,G)$ is well defined and $C^M$.
\end{enumerate}
Analogous results hold for spaces of multilinear mappings. \end{proposition}

\demo{Proof}
\thetag{1}
That the first three topologies on $L(E,F)$ have the same bounded sets has been shown
in \cite[5.3, 5.18]{KM97}.
The inclusion $C^M(E,F)\to C^\infty(E,F)$ is bounded by the uniform boundedness principle \cite[5.18]{KM97}.
Conversely, the inclusion $L(E,F)\to C^M(E,F)$ is bounded
by the uniform boundedness principle \thetag{\ref{nmb:2.9}}.

\thetag{2} The assertion for $C^\infty$ is true by \cite[3.12]{KM97}
since $L(E,F)$ is closed in $C^\infty(E,F)$.

If $f$ is $C^M$ let $L\geq M$ be a non-quasianalytic weight-sequence
and let $c:\mathbb R\to U$ be a $C^L$-curve. We have to show that $f^\vee\circ c$ is 
$C^L$ into $L(F,G)$. By the uniform boundedness principle \cite[5.18]{KM97} and
\thetag{\ref{nmb:2.3}} it suffices to show that 
$t\mapsto \ell\bigl(f^\vee(c(t))(v)\bigr)=\ell\bigl(f(c(t),v)\bigr)\in \mathbb R$ is $C^L$ for each
$\ell\in G^*$ and $v\in F$; this is obviously true. 
Conversely, let $f^\vee:U\to L(F,G)$ be $C^M$ and let $L\geq M$ be a non-quasianalytic weight-sequence. 
We claim that $f:U\times F\to G$ is $C^L$. By composing with $\ell\in G^*$ we
may assume that $G=\mathbb R$.
By induction we have \begin{align*}
&d^kf(x,w_0)\big((v_k,w_k),\dots,(v_1,w_1)\big) = d^k (f^\vee)(x)(v_k,\dots,v_1)(w_0) +\\&
+ \sum_{i=1}^k d^{k-1}(f^\vee)(x)(v_k,\dots,\widehat{v_i,}\dots,v_1)(w_i)
\end{align*}
We check condition \thetag{\ref{nmb:2.7}.4} for $f$ where $x\in K$ which is compact in $U$:
\begin{align*}
&\|d^k f (x,w_0)\|_{L^k(E_B\times F_{B'},\mathbb R)}\le
\\&
\le \|d^k (f^\vee)(x)(\dots)(w_0)\|_{L^k(E_B,\mathbb R)}
+ \sum_{i=1}^k \|d^{k-1}(f^\vee)(x)\|_{L^{k-1}(E_B,L(F_{B'},\mathbb R))}
\\&
\le \|d^k (f^\vee)(x)\|_{L^k(E_B,L(F_{B'},\mathbb R))}\|w_0\|_{B'}
+ \sum_{i=1}^k \|d^{k-1}(f^\vee)(x)\|_{L^{k-1}(E_B,L(F_{B'},\mathbb R))}
\\&
\le C\, \rho^k\, k!\, L_k \|w_0\|_{B'}
+ \sum_{i=1}^k C \,\rho^{k-1}\, (k-1)!\, L_{k-1} = C\,\rho^k \,k!\,L_k \Bigl(\|w_0\|_{B'}+ \tfrac{L_{k-1}}{\rho\, L_k}\Bigr)
\end{align*}
where we used \thetag{\ref{nmb:2.7}.4} for $L(i_{B'},\mathbb R)\circ f^\vee: U\to
L(F_{B'},\mathbb R)$.
Since $L$ is increasing, $f$ is $C^L$.
\qed\enddemo

\begin{theorem}\label{nmb:4.3}
Let $Q=(Q_k)$ be an $\mathcal L$-intersectable quasianalytic weight sequence.
Let $U\subseteq E$ be $c^\infty$-open in a convenient vector space, let $F$
be another convenient vector space, %
and $f:U\to F$ a mapping. Then the following statements are equivalent:
\begin{enumerate}
\item[(1)] $f$ is $C^Q$, i.e., for all $L\in\mathcal L(Q)$ we have $f\circ c$ is $C^L$ for all $C^L$-curves $c$.
\item[(2)] $f|_{U\cap E_B}:E_B\supseteq U\cap E_B\to F$ is $C^Q$ for each closed bounded absolutely convex $B$ in $E$.
\item[(3)] For all $L\in\mathcal L(Q)$ the curve $f\circ c$ is $C^L$ for all $C_b^L$-curves $c$.
\item[(4)] $\pi_V\circ f$ is $C_b^Q$ for all absolutely convex 0-neighborhoods $V$ in $F$
and the associated mapping $\pi_V:F\to F_V$.
\end{enumerate}
\end{theorem}
\begin{demo}{Proof}
This follows from \thetag{\ref{nmb:2.8}} for $\overline{L}:=L$ since $C^Q:=\bigcap_{L\in\mathcal{L}(Q)}C^L$
and $C_b^Q=\bigcap_{L\in\mathcal{L}(Q)}C_b^L$.%
\qed\end{demo}

\begin{theorem}[cf.\ {\cite[4.4]{KMRc}}]\label{nmb:4.4}
Let $Q=(Q_k)$ be an $\mathcal L$-intersectable quasianalytic weight sequence.
Let $E$ and $F$ be convenient vector spaces and let $U\subseteq E$ be $c^\infty$-open.
Then
\[
C^Q(U,F)
=\varprojlim_{L\in\mathcal{L}(Q),c\in C^L} C^L(\mathbb R,F)
=\varprojlim_{B\subseteq E} C^Q(U\cap E_B,F)
=\varprojlim_{L\in\mathcal{L}(Q),s\in C_b^L} C^L(\mathbb R,F)
\]
as vector spaces with bornology, where
$c$ runs through all $C^L$-curves in $U$ for $L\in\mathcal{L}(Q)$,
$B$ runs through all bounded closed absolutely convex subsets of $E$,
and $s$ runs through all $C_b^L$-curves in $U$ for $L\in\mathcal{L}(Q)$.
\end{theorem}
\begin{demo}{Proof}
This follows by applying $\varprojlim_{L\in\mathcal{L}(Q)}$ to \thetag{\ref{nmb:2.11}}.
\qed\end{demo}

\subsection{\label{nmb:4.5} Jet spaces}
Let $E$ and $F$ be Banach spaces
and $A\subseteq E$ convex.
We consider the linear space
$C^\infty(A,F)$ consisting of all sequences $(f^k)_{k}\in\prod_{k\in\mathbb N} C(A,L^k(E,F))$ satisfying
\[f^k(y)(v)-f^k(x)(v) =
\int_0^1 f^{k+1}(x+t(y-x))(y-x,v)\,dt \]
for all $k\in\mathbb N$, $x,y\in A$, and $v\in E^k$.
If $A$ is open we can identify this space with that of all smooth functions
$A\to F$ by the passage to jets.

In addition, let $M=(M_k)$ be a weight sequence and $(r_k)$ a sequence of positive real numbers. 
Then we consider the normed spaces \[
C^M_{(r_k)}(A,F):=\Bigl\{(f^k)_k\in C^\infty(A,F):\|(f^k)\|_{(r_k)}<\infty\Bigr\}
\]
where the norm is given by
\[
\|(f^k)\|_{(r_k)}:=\sup\Bigl\{\frac{\|f^k(a)(v_1,\dots,v_k)\|}{k!\,r_k\,M_k\,\|v_1\|\cdot\dots\cdot \|v_k\|}:
k\in\mathbb N, a\in A, v_i\in E\Bigr\}.
\]
If $(r_k)=(\rho^k)$ for some $\rho>0$ we just write $\rho$ instead of $(r_k)$ as indices.
The spaces $C^M_{(r_k)}(A,F)$ are Banach spaces, since they are closed
in $\ell^\infty(\mathbb N,\ell^\infty(A,L^k(E,F)))$ via $(f^k)_k\mapsto (k\mapsto \frac1{k!\,r_k\,M_k} f^k)$.

If $A$ is open, $C^\infty(A,F)$ and $C_\rho^M(A,F)$ coincide with the convenient spaces treated before.

\begin{theorem}[cf.\ {\cite[4.6]{KMRc}}]\label{nmb:4.6} Let $M=(M_k)$ be a non-quasianalytic weight sequence
or an $\mathcal L$-intersectable quasianalytic weight sequence.
Let $E$ and $F$ be Banach spaces and let $U\subseteq E$ be open and convex.
Then the space $C^M(U,F)=C^M_b(U,F)$ can be described bornologically in the following equivalent 
ways, i.e., these constructions give the same
vector space and the same bounded sets
\begin{align*}
\varprojlim_K \varinjlim_{\rho,W} C^M_\rho(W,F) \tag{1} \\
\varprojlim_K \varinjlim_\rho C^M_\rho(K,F) \tag{2} \\
\varprojlim_{K,(r_k)} C^M_{(r_k)}(K,F) \tag{3} \\
\end{align*}
Moreover, all involved inductive limits are regular, i.e.\ the bounded sets of 
the inductive limits are contained and bounded in some step.

Here $K$ runs through all compact convex subsets of $U$ ordered by inclusion, 
$W$ runs through the open subsets $K\subseteq W\subseteq U$ again ordered by inclusion,
$\rho$ runs through the positive real numbers,
$(r_k)$ runs through all sequences of positive real numbers for which $\rho^k/r_k\to 0$ for all $\rho>0$.
\end{theorem}

\begin{demo}{Proof} This proof is almost identical with that of \cite[4.6]{KMRc}.
The only change is to use \thetag{\ref{nmb:2.7}} and \thetag{\ref{nmb:4.3}} 
instead of \cite[3.9]{KMRc} to show that all these descriptions give $C^M(U,F)$ as vector space.
\qed\end{demo}

\begin{lemma}[cf.\ {\cite[4.7]{KMRc}}]\label{nmb:4.7}
Let $M$ be a non-quasianalytic weight sequence. For any convenient
vector  space  $E$ the flip of variables induces an isomorphism
$L(E,C^M(\mathbb  R,\mathbb  R))  \cong  C^M(\mathbb  R,E')$ as
vector spaces.
\end{lemma}

\begin{demo}{Proof}  This  proof  is  identical  with  that  of
\cite[4.7]{KMRc}  but  uses  \thetag{\ref{nmb:2.9}}  instead of
\cite[4.1]{KMRc}    and   \thetag{\ref{nmb:2.3}}   instead   of
\cite[3.5]{KMRc}. %
\qed\end{demo}

\begin{lemma}[cf.\ {\cite[4.8]{KMRc}}]\label{nmb:4.8}
Let $M=(M_k)$ be a non-quasianalytic weight sequence.
By  $\lambda^M(\mathbb  R)$ we denote the $c^\infty$-closure of
the linear subspace
generated  by  
$\{\on{ev}_t:t\in\mathbb  R\}$  in  $C^M(\mathbb R,\mathbb  R)'$  and  let 
$\delta:\mathbb R\to\lambda^M(\mathbb R)$
be given by $t\mapsto \on{ev}_t$.
Then $\lambda^M(\mathbb R)$ is the free convenient vector space over $C^M$, i.e.\
for every convenient vector space $G$ the $C^M$-curve $\delta$ induces a
bornological isomorphism
\begin{align*}
\delta^*:L(\lambda^M(\mathbb R),G)&\cong C^M(\mathbb R,G).
\end{align*}
\end{lemma}
We expect $\lambda^M(\mathbb R)$ to be equal to $C^M(\mathbb R,\mathbb R)'$ as it is the case for the analogous situation of smooth mappings, see \cite[23.11]{KM97}, and of
holomorphic mappings, see \cite{Siegl95} and \cite{Siegl97}.

\begin{demo}{Proof} The proof goes along the same lines as in \cite[23.6]{KM97} and in \cite[5.1.1]{FK88}.
It is identical with that of \cite[4.8]{KMRc} but uses \thetag{\ref{nmb:2.3}}, \thetag{\ref{nmb:2.9}}, 
and \thetag{\ref{nmb:4.2}} in that order.
\qed\end{demo}

\begin{corollary}[cf.\ {\cite[4.9]{KMRc}}]\label{nmb:4.9}
Let $L=(L_k)$ and $L'=(L'_k)$ be non-quasianalytic weight sequences.
We have the following isomorphisms of linear spaces
\begin{enumerate}
\item[(1)] $C^\infty(\mathbb R,C^L(\mathbb R,\mathbb R)) \cong C^L(\mathbb R,C^\infty(\mathbb R,\mathbb R))$ 
\item[(2)] $C^\omega(\mathbb R,C^L(\mathbb R,\mathbb R)) \cong C^L(\mathbb R,C^\omega(\mathbb R,\mathbb R))$
\item[(3)] $C^{L'}(\mathbb R,C^L(\mathbb R,\mathbb R)) \cong C^L(\mathbb R,C^{L'}(\mathbb R,\mathbb R))$
\end{enumerate}
\end{corollary}

\begin{demo}{Proof} This proof is that of \cite[4.9]{KMRc} with other refernces:
For $\alpha\in\{\infty,\omega,L'\}$ we get
\begin{align*}
C^L(\mathbb R,C^\alpha(\mathbb R,\mathbb R)) &\cong
L(\lambda^L(\mathbb R),C^\alpha(\mathbb R,\mathbb R))\qquad\text{ by \thetag{\ref{nmb:4.8}}}\\
&\cong C^\alpha(\mathbb R,L(\lambda^L(\mathbb R),\mathbb R))\qquad\text{by \thetag{\ref{nmb:4.7}}, 
\cite[3.13.4, 5.3, 11.15]{KM97}}\\
&\cong C^\alpha(\mathbb R,C^L(\mathbb R,\mathbb R))\qquad\text{ by \thetag{\ref{nmb:4.8}}.}\qed
\end{align*}
\end{demo}

\begin{theorem}[Canonical isomorphisms]\label{nmb:4.10}
Let $M=(M_k)$ be a non-quasianalytic weight sequences or an 
$\mathcal L$-intersectable quasianalytic weight-sequences; 
likewise $M'=(M'_k)$. Let $E$, $F$ be convenient vector spaces and let
$W_i$ be $c^\infty$-open subsets in such. We have the following natural bornological isomorphisms:
\begin{enumerate}
\item[(1)] $C^M(W_1,C^{M'}(W_2,F))\cong C^{M'}(W_2,C^M(W_1,F))$,
\item[(2)] $C^M(W_1,C^\infty(W_2,F))\cong C^\infty(W_2,C^M(W_1,F))$.
\item[(3)] $C^M(W_1,C^\omega(W_2,F))\cong C^\omega(W_2,C^M(W_1,F))$.
\item[(4)] $C^M(W_1,L(E,F))\cong L(E,C^M(W_1,F))$.
\item[(5)] $C^M(W_1,\ell^\infty(X,F))\cong \ell^\infty(X,C^M(W_1,F))$.
\item[(6)] $C^M(W_1,\operatorname{\mathcal Lip}^k(X,F))\cong \operatorname{\mathcal Lip}^k(X,C^M(W_1,F))$.
\end{enumerate}
In \thetag{5} the space $X$ is an $\ell^\infty$-space, i.e.\ 
a set together with a bornology induced by a family of real valued functions on $X$, cf.\ \cite[1.2.4]{FK88}.
In \thetag{6} the space $X$ is a $\operatorname{\mathcal Lip}^k$-space, cf.\ 
\cite[1.4.1]{FK88}.
The spaces $\ell^\infty(X,F)$ and $\operatorname{\mathcal Lip}^k(W,F)$ are defined in \cite[3.6.1 and 4.4.1]{FK88}. 
\end{theorem}
\begin{demo}{Proof} This proof is very similar with that of \cite[4.8]{KMRc} but written differently.
Let $\mathcal C^1$ and $\mathcal C^2$ denote any of the functions spaces mentioned above
and $X_1$ and $X_2$ the corresponding domains.
In order to show that the flip of coordinates $f\mapsto \tilde f$, 
$\mathcal C^1(X_1,\mathcal C^2(X_2,F))\to \mathcal C^2(X_2,\mathcal C^1(X_1,F))$
is a well-defined bounded linear mapping we have to show:
\begin{itemize}
\item $\tilde f(x_2)\in\mathcal C^1(X_1,F)$, which is obvious, since
$\tilde f(x_2)=\on{ev}_{x_2}\circ f:X_1\to \mathcal C^2(X_2,F)\to F$.
\item $\tilde f\in\mathcal C^2(X_2,\mathcal C^1(X_1,F))$, which we will show below.
\item $f\mapsto \tilde f$ is bounded and linear, which follows by applying the appropriate
uniform boundedness theorem for $\mathcal C^2$ and $\mathcal C^1$ since
$f\mapsto \on{ev}_{x_1}\circ \on{ev}_{x_2}\circ\tilde f=\on{ev}_{x_2}\circ\on{ev}_{x_1}\circ f$ is bounded and linear.
\end{itemize}
All occurring function spaces are convenient and satisfy the uniform $\mathcal{S}$-boundedness theorem, 
where $\mathcal{S}$ is the set of point evaluations:
\begin{itemize}
\item[$C^M$] by \thetag{\ref{nmb:1.10}} and \thetag{\ref{nmb:2.9}}.
\item[$C^\infty$] by \cite[2.14.3, 5.26]{KM97}, \item[$C^\omega$] by \cite[11.11, 11.12]{KM97}, 
\item[$L$] by \cite[2.14.3, 5.18]{KM97}, \item[$\ell^\infty$] by \cite[2.15, 5.24, 5.25]{KM97} or 
\cite[3.6.1 and 3.6.6]{FK88}
\item[$\operatorname{\mathcal Lip}^k$] by \cite[4.4.2 and 4.4.7]{FK88}
\end{itemize}

It remains to check that $\tilde f$ is of the appropriate class:
\begin{itemize}
\item[\thetag{1}] follows by composing with the appropriate (non-quasianalytic) curves
$c_1:\mathbb R\to W_1$, $c_2:\mathbb R\to W_2$ and $\lambda\in F^*$
and thereby reducing the statement to the special case in \thetag{\ref{nmb:4.9}.3}.
\item[\thetag{2}] as for \thetag{1} using \thetag{\ref{nmb:4.9}.1}.
\item[\thetag{3}] follows
by composing with $c_2\in C^{\beta_2}(\mathbb R,W_2)$, where $\beta_2$ is in 
$\{\infty,\omega\}$, and with $C^L(c_1,\lambda):C^M(W_1,F) \to C^L(\mathbb R,\mathbb R)$
where $c_1 \in C^L(\mathbb R,W_1)$ with $L\geq M$ non-quasianalytic and
$\lambda\in F^*$.
Then $C^L(c_1,\lambda)\circ \tilde f \circ c_2 =
(C^{\beta_2}(c_2,\lambda)\circ f \circ c_1)^\sim
: \mathbb R \to C^L(\mathbb R,\mathbb R)$
is $C^{\beta_2}$ by \thetag{\ref{nmb:4.9}.1} and \thetag{\ref{nmb:4.9}.2}, since $C^{\beta_2}(c_2,\lambda)\circ f \circ c_1
: \mathbb R \to W_1 \to C^\omega(W_2,F) \to C^{\beta_2}(\mathbb R,\mathbb R)$
is $C^L$.
\newline
For the inverse, compose with $c_1$ and $C^{\beta_2}(c_2,\lambda):C^\omega(W_2,F) \to C^{\beta_2}(\mathbb R,\mathbb R)$.
Then $C^{\beta_2}(c_2,\lambda)\circ \tilde f \circ c_1 =
(C^L(c_1,\lambda)\circ f \circ c_2)^\sim
: \mathbb R \to C^{\beta_2}(\mathbb R,\mathbb R)$
is $C^L$ by \thetag{\ref{nmb:4.9}.1} and \thetag{\ref{nmb:4.9}.2}, since $C^L(c_1,\lambda)\circ f \circ c_2
: \mathbb R \to W_2 \to C^L(W_1,F) \to C^L(\mathbb R,\mathbb R)$
is $C^{\beta_2}$.
\item[\thetag{4}] since $L(E,F)$ is the $c^\infty$-closed subspace of $C^M(E,F)$ formed
by the linear $C^M$-mappings.
\item[\thetag{5}] follows from
\thetag{4}, using the free convenient vector spaces $\ell^1(X)$ over the $\ell^\infty$-space $X$, see \cite[5.1.24 or 5.2.3]{FK88}, satisfying $\ell^\infty(X,F)\cong L(\ell^1(X),F)$.
\item[\thetag{6}]
follows from
\thetag{4}, using the free convenient vector spaces $\lambda^k(X)$ over the 
$\operatorname{\mathcal Lip}^k$-space $X$, satisfying 
$\operatorname{\mathcal Lip}^k(X,F)\cong L(\lambda^k(X),F)$. 
Existence of this free convenient vector space can be proved in a similar way as in \thetag{\ref{nmb:4.8}}.
\qed
\end{itemize}
\end{demo}

\section{\label{nmb:5} Manifolds of quasianalytic mappings}

For  manifolds  of real analytic mappings \cite{KrieglMichor90}
we  could  prove  that  composition and inversion (on groups of
real  analytic diffeomorphisms) are again $C^\omega$ by testing
along  $C^\infty$-curves and $C^\omega$-curves separately. Here
this  does  not  (yet) work. We have to test along $C^L$-curves
for  all  $L$  in  $\mathcal L(Q)$, but for those $L$ we do not
have  cartesian  closedness in general. But it suffices to test
along  $C^Q$-mappings from open sets in Banach spaces, and this
is a workable replacement.

\subsection{\label{nmb:5.1}$C^Q$-manifolds}
Let $Q=(Q_k)$ be an $\mathcal L$-intersectable quasianalytic weight sequence of moderate growth.
A $C^Q$-manifold is a smooth manifold such that all chart changings are
$C^Q$-mappings. Likewise for $C^Q$-bundles and $C^Q$ Lie groups. 
Note that any finite dimensional (always assumed paracompact) $C^\infty$-manifold admits a
$C^\infty$-diffeomorphic real analytic structure thus also a
$C^Q$-structure. Maybe, any finite dimensional $C^Q$-manifold admits a
$C^Q$-diffeomorphic real analytic structure.
This would follow from:

\subsection*{Conjecture} {\it
Let  $X$  be  a  finite  dimensional  real  analytic  manifold.
Consider the space $C^Q(X,\mathbb R)$ of all $C^Q$-functions on
$X$,  equipped  with the (obvious) Whitney $C^Q$-topology. Then
$C^\omega(X,\mathbb R)$ is dense in $C^Q(X,\mathbb R)$. }

This conjecture is the analogon of \cite[Proposition 9]{Grauert58}. 
\subsection{\label{nmb:5.2} Banach plots}
Let $Q=(Q_k)$ be an $\mathcal L$-intersectable quasianalytic weight sequence of moderate growth.
Let  $X$  be  a $C^Q$-manifold. By a {\it $C^Q$-plot in $X$} we
mean  a $C^Q$-mapping $c:D\to X$ where $D\subset E$ is the open
unit ball in a Banach space $E$.

\begin{lemma*}
A mapping between $C^Q$-manifolds is $C^Q$ if and only if it maps $C^Q$-plots to $C^Q$-plots. \end{lemma*}

\begin{demo}{Proof}
For  a  convenient  vector space $E$ the $c^\infty$-topology is
the final topology for all injections $E_B\to E$ where $B$ runs
through  all  closed  absolutely convex bounded subsets of $E$.
The $c^\infty$-topology on a $c^\infty$-open subset $U\subseteq
E$  is  final  with respect to all injections $E_B\cap U\to U$.
For  a  $C^Q$-manifold  the  topology  is the final one for all
$C^Q$-plots.  Let  $f:X\to  Y$  be the mapping. If $f$ respects
$C^Q$-plots  it  is continuous and so we may assume that $Y$ is
$c^\infty$-open  in  a  convenient  vector  space  $F$ and then
likewise  for  $X\subseteq  E$.  The  (affine) plots induced by
$X\cap E_B\subset X$ are $C^Q$.
By  definition  $f$ is $C^Q$ if and only if it is $C^L$ for all
$L\in\mathcal  L(Q)$  and  this  is the case if $f$ is $C^L$ on
$X\cap    E_B$   for   all   $B$   by   \thetag{\ref{nmb:2.8}}.
\qed\end{demo}

\subsection{\label{nmb:5.3}Spaces of $C^Q$-sections}
Let $p:E\to B$ be a $C^Q$ vector bundle (possibly infinite dimensional).
The space $C^Q(B\gets E)$ of all $C^Q$-sections is a convenient
vector space with the structure induced by \begin{gather*}
C^Q(B\gets E) \to \prod_\alpha C^Q(u_\alpha(U_\alpha),V)
\\
s\mapsto \on{pr}_2\circ \psi_\alpha \circ s \circ u_\alpha^{-1}
\end{gather*}
where 
$B\supseteq U_\alpha -\raisebox{0.1pt}{$\mkern-16mu\frac{\;\;u_\alpha\;}{\;\;\;}\mkern-16mu$}
\to u_\alpha(U_\alpha)\subseteq W$ is a
$C^Q$-atlas for $B$ which we assume to be modeled on a convenient vector
space $W$, and where
$\psi_\alpha :E|_{U_\alpha}\to U_\alpha\times V$ form a vector bundle atlas over
charts $U_\alpha$ of $B$. 

\begin{lemma*}
Let  $D$  be  a  unit ball in a Banach space. A mapping $c:D\to
C^Q(B\gets  E)$  is  a  $C^Q$-plot  if  and  only if $c^\wedge:
D\times B\to E$ is $C^Q$.
\end{lemma*}

\begin{demo}{Proof}
By the description of the structure on $C^Q(B\gets E)$ we may assume that
$B$  is  $c^\infty$-open  in  a convenient vector space $W$ and
that  $E=B\times  V$.  Then we have $C^Q(B\gets B\times V)\cong
C^Q(B,V)$. Thus the statement follows from the
exponential law \thetag{\ref{nmb:3.3}}. \qed\end{demo}

Let  $U\subseteq  E$  be  an  open neighborhood of $s(B)$ for a
section  $s$  and  let $q:F\to B$ be another vector bundle. The
set  $C^Q(B\gets  U)$  of  all  $C^Q$-sections $s':B\to E$ with
$s'(B)\subset  U$  is  open  in  the  convenient  vector  space
$C^Q(B\gets E)$ if $B$ is compact.
An immediate consequence of the lemma is the following: If $U\subseteq E$ is an open
neighborhood of $s(B)$ for a section $s$, $F\to B$ is another vector bundle and if $f:U\to
F$ is a fiber respecting $C^Q$-mapping, then $f_*:C^Q(B\gets U)\to
C^Q(B\gets F)$ is $C^Q$ on the open neighborhood $C^Q(B\gets U)$ of $s$ in
$C^Q(B\gets E)$. We have $(d(f_*)(s)v)_x = d(f|_{U\cap E_x})(s(x))(v(x))$.

\begin{theorem}\label{nmb:5.4}
Let  $Q=(Q_k)$  be  an $\mathcal L$-intersectable quasianalytic
weight sequence of moderate growth.
Let  $A$ and $B$ be finite dimensional $C^Q$-manifolds with $A$
compact  and $B$ equipped with a $C^Q$ Riemann metric. Then the
space $C^Q(A,B)$ of all $C^Q$-mappings $A\to B$ is a
$C^Q$-manifold  modeled on convenient vector spaces $C^Q(A\gets
f^*TB)$ of
$C^Q$-sections  of pullback bundles along $f:A\to B$. Moreover,
a  mapping  $c:D\to  C^Q(A,B)$  is  a $C^Q$-plot if and only if
$c^\wedge :D\times A\to B$ is $C^Q$.
\end{theorem}

If  the  $C^Q$-structure  on  $B$ is induced by a real analytic
structure  then  there  exists  a  real analytic Riemann metric
which in turn is $C^Q$.

\begin{demo}{Proof}
$C^Q$-vector fields have $C^Q$-flows by
\cite{Komatsu80};  applying  this  to the geodesic spray we get
the $C^Q$ exponential mapping $\exp: TB\supseteq U\to B$ of the
Riemann metric,
defined on a suitable open neighborhood of the zero section. We
may assume that $U$ is chosen in such a way that
$(\pi_B,\exp):U\to B\times B$ is a $C^Q$-diffeomorphism onto
an  open neighborhood $V$ of the diagonal, by the $C^Q$ inverse
function theorem
due to \cite{Komatsu79}. 
For $f\in C^Q(A,B)$ we consider the pullback vector bundle
\[\xymatrix{
A\times TB & A\times_BTB \ar@{_(->}[l] \ar@{=}[r] & f^*TB \ar[r]^{\pi_B^*f} \ar[d]_{f^*\pi_B} & TB \ar[d]^{\pi_B} \\
& & A \ar[r]^f & B
}\]
Then  the  convenient  space of sections $C^Q(A\gets f^*TB)$ is
canonically  isomorphic  to  the  space  $C^Q(A,TB)_f:=  \{h\in
C^Q(A,TB):\pi_B\circ  h=f\}$  via  $s\mapsto (\pi_B^*f)\circ s$
and $(\on{Id}_A,h)\mapsfrom h$. %
Now let 
\begin{gather*} 
U_f :=\{g\in C^Q(A,B):(f(x),\; g(x))\in V
\text{ for all }x\in A\},\\
u_f:U_f\to C^Q(A\gets f^*TB),\\
u_f(g)(x) = (x,\exp_{f(x)}^{-1}(g(x))) = (x,((\pi_B,\exp)^{-1}\circ(f,g))(x)).
\end{gather*}
Then  $u_f:U_f  \to  \{s\in  C^Q(A\gets  f^*TB):  s(A)\subseteq
f^*U=(\pi_B^*f)^{-1}(U)\}$   is   a   bijection   with  inverse
$u_f^{-1}(s)  =  \exp\circ(\pi_B^*f)\circ  s$, where we view $U
\to B$ as a fiber bundle. The set $u_f(U_f)$ is open in
$C^Q(A\gets f^*TB)$ for the topology described above in \thetag{\ref{nmb:5.3}} since $A$ is compact and
the push forward $u_f$ is $C^Q$ since it respects $C^Q$-plots
by lemma \thetag{\ref{nmb:5.3}}. 

Now we consider the atlas $(U_f,u_f)_{f\in C^Q(A,B)}$ for
$C^Q(A,B)$. Its chart change mappings are given for 
$s\in u_g(U_f\cap U_g)\subseteq C^Q(A\gets g^*TB)$ by 
\begin{align*} (u_f\circ u_g^{-1})(s) &= (\on{Id}_A,(\pi_B,\exp)^{-1}\circ(f,\exp\circ(\pi_B^*g)\circ s)) \\
&= (\tau_f^{-1}\circ\tau_g)_*(s),
\end{align*}
where $\tau_g(x,Y_{g(x)}) := (x,\exp_{g(x)}(Y_{g(x)}))$
is a $C^Q$-diffeomorphism $\tau_g:g^*TB \supseteq g^*U \to (g\times \on{Id}_B)^{-1}(V)\subseteq A\times B$
which is fiber respecting over $A$. 
The chart change $u_f\circ u_g^{-1} = (\tau_f^{-1}\circ \tau_g)_*$ is defined on an open
subset and it is also $C^Q$ since it respects $C^Q$-plots
by lemma \thetag{\ref{nmb:5.3}}. 
Finally   for   the   topology   on   $C^Q(A,B)$  we  take  the
identification    topology    from   this   atlas   (with   the
$c^\infty$-topologies   on   the  modeling  spaces),  which  is
obviously finer than the
compact-open topology and thus Hausdorff.

The equation $u_f\circ u_g^{-1} = (\tau_f^{-1}\circ \tau_g)_*$ shows that
the $C^Q$-structure does not depend on the choice of the
$C^Q$ Riemannian metric on $B$.

The statement on $C^Q$-plots follows from lemma \thetag{\ref{nmb:5.3}}. \qed\end{demo}

\begin{corollary}\label{nmb:5.5}
Let $A_1,A_2$ and $B$ be finite dimensional $C^Q$-manifolds with $A_1$ and
$A_2$ compact. Then composition 
\[
C^Q(A_2,B) \times C^Q(A_1,A_2) \to C^Q(A_1,B), \quad (f,g) \mapsto f\circ g
\]
is  $C^Q$.  However,  if  $N=(N_k)$  is another weight sequence
($\mathcal L$-intersectable quasianalytic)
with $(N_k/Q_k)^{1/k} \searrow 0$ then composition is {\bf not} $C^N$.
\end{corollary}

\demo{Proof}
Composition maps $C^Q$-plots to $C^Q$-plots, so it is $C^Q$. 
Let  $A_1=A_2=S^1$  and  $B=\mathbb  R$.  Then by \cite[Theorem
1]{Thilliez08} or \cite[2.1.5]{KMRc} there exists $f\in
C^Q(S^1,\mathbb  R)\setminus  C^N(S^1,\mathbb  R)$. We consider
$f$ as a periodic function $\mathbb R\to \mathbb R$.
The  universal covering space of $C^Q(S^1,S^1)$ consists of all
$2\pi\mathbb  Z$-equivariant mappings in $C^Q(\mathbb R,\mathbb
R)$,  namely  the  space  of  all  $g+\on{Id}_{\mathbb  R}$ for
$2\pi$-periodic  $g\in  C^Q$.  Thus  $C^Q(S^1,S^1)$  is  a real
analytic  manifold and $t\mapsto (x\mapsto x+t)$ induces a real
analytic  curve  $c$ in $C^Q(S^1,S^1)$. But $f_*\circ c$ is not
$C^N$ since:
\begin{align*}
\frac{(\partial_t^k|_{t=0}(f_*\circ c)(t))(x)}{k!\rho^k N_k} =
\frac{\partial_t^k|_{t=0} f(x+t)}{k!\rho^k N_k} = \frac{f^{(k)}(x)}{k!\rho^k N_k}
\end{align*}
which is unbounded in $k$ for $x$ in a suitable compact set and
for all $\rho>0$, since $f\notin C^N$.
\qed\enddemo

\begin{theorem}\label{nmb:5.6}
Let $Q=(Q_k)$ be a, $\mathcal L$-intersectable quasianalytic weight sequence of moderate growth.
Let $A$ be a compact ($\Rightarrow$ finite dimensional) $C^Q$-manifold.
Then the group $\on{Diff}^Q(A)$ of all $C^Q$-diffeomorphisms of $A$ is an
open subset of the $C^Q$-manifold $C^Q(A,A)$. Moreover, it is
a $C^Q$-regular $C^Q$ Lie group: Inversion and composition are $C^Q$. Its Lie algebra
consists of all $C^Q$-vector fields on $A$, with the negative of the usual
bracket as Lie bracket. The exponential mapping is $C^Q$. It is not surjective onto any neighborhood
of $\on{Id}_A$. \end{theorem}

Following \cite{KM97r}, see also \cite[38.4]{KM97}, a $C^Q$-Lie
group  $G$  with  Lie  algebra  $\mathfrak  g=T_eG$  is  called
$C^Q$-regular if the following holds:
\begin{itemize}
\item  For  each  $C^Q$-curve $X\in C^Q(\mathbb R,\mathfrak g)$
there  exists a $C^Q$-curve $g\in C^Q(\mathbb R,G)$ whose right
logarithmic derivative is $X$, i.e.,
\[
\begin{cases} g(0) &= e \\
\partial_t g(t) &= T_e(\mu^{g(t)})X(t) = X(t).g(t)
\end{cases} \]
The curve $g$ is uniquely determined by its initial value $g(0)$, if it
exists.
\item
Put  $\on{evol}^r_G(X)=g(1)$  where  $g$ is the unique solution
required  above.  Then  $\on{evol}^r_G: C^Q(\mathbb R,\mathfrak
g)\to G$ is required to be
$C^Q$ also. \end{itemize}

\demo{Proof}
The group $\on{Diff}^Q(A)$ is open in $C^Q(A,A)$ since it is open in the
coarser  $C^1$ compact-open topology, see \cite[43.1]{KM97}. So
$\on{Diff}^Q(A)$  is  a $C^Q$-manifold and composition is $C^Q$
by
\thetag{\ref{nmb:5.4}} and \thetag{\ref{nmb:5.5}}. To show that inversion is $C^Q$ let $c$ be a
$C^Q$-plot in $\on{Diff}^Q(A)$. By \thetag{\ref{nmb:5.4}} the map $c^\wedge: D\times A\to A$ is
$C^Q$ and $(\on{inv}\circ\, c)^\wedge:D\times A\to A$ satisfies the Banach manifold
implicit equation $c^\wedge(t,(\on{inv}\circ\, c)^\wedge(t,x))=x$ for $x\in A$. By the Banach
$C^Q$  implicit  function theorem \cite{Yamanaka89} the mapping
$(\on{inv}\circ\,  c)^\wedge$  is locally $C^Q$ and thus $C^Q$.
By \thetag{\ref{nmb:5.4}}
again, $\on{inv}\circ\, c$ is a $C^Q$-plot in $\on{Diff}^Q(A)$.
So  $\on{inv}:\on{Diff}^Q(A)\to  \on{Diff}^Q(A)$  is $C^Q$. The
Lie  algebra of $\on{Diff}^Q(A)$ is the convenient vector space
of  all  $C^Q$-vector  fields  on $A$, with the negative of the
usual Lie bracket
(compare with the proof of \cite[43.1]{KM97}). 
To show that $\on{Diff}^Q(A)$ is a $C^Q$-regular Lie group, we choose a
$C^Q$-plot in the space of $C^Q$-curves in the Lie algebra of all $C^Q$
vector fields on $A$, $c:D\to C^Q(\mathbb R,C^Q(A\gets TA))$. By
lemma \thetag{\ref{nmb:5.3}} $c$ corresponds to a $(D\times \mathbb R)$-time-dependent $C^Q$
vector field $c^{\wedge \wedge }:D\times\mathbb R\times A\to TA$. Since $C^Q$-vector
fields   have   $C^Q$-flows   and   since   $A$   is   compact,
$\on{evol}^r(c^\wedge  (s))(t)  =  \on{Fl}^{c^\wedge (s)}_t$ is
$C^Q$ in all
variables by \cite{Yamanaka91}.
Thus $\on{Diff}^Q(A)$ is a $C^Q$-regular $C^Q$ Lie group.

The exponential mapping is $\on{evol}^r$ applied to constant curves in the
Lie  algebra,  i.e.,  it  consists of flows of autonomous $C^Q$
vector  fields. That the exponential map is not surjective onto
any $C^Q$-neighborhood of
the identity follows from \cite[43.5]{KM97} for $A=S^1$. This example can
be embedded into any compact manifold, see \cite{Grabowski88}. \qed\enddemo

\def\cprime{$'$}
\providecommand{\bysame}{\leavevmode\hbox to3em{\hrulefill}\thinspace}
\providecommand{\MR}{\relax\ifhmode\unskip\space\fi MR }
\providecommand{\MRhref}[2]{%
  \href{http://www.ams.org/mathscinet-getitem?mr=#1}{#2}
}
\providecommand{\href}[2]{#2}

\end{document}